\documentclass[a4paper,10pt,leqno]{amsart}

\usepackage[latin1]{inputenc}
\usepackage[english]{babel}
\usepackage{amsthm,amssymb,amsmath,amsopn,amsfonts,pb-diagram,footmisc,amsthm,amsaddr,mathrsfs,mathtools}
\usepackage{stmaryrd,calc,enumerate,mathrsfs,amscd,footnote,hyperref}
\usepackage{graphicx,color}
\usepackage{fancyhdr}	
\usepackage{indentfirst}
\usepackage{geometry}
\hypersetup{
	colorlinks=true,
	linkcolor=blue,
	citecolor=red,
	urlcolor=blue,
	linktoc=all}
	
\DeclareMathOperator{\dist}{dist}

\DeclareMathOperator{\pv}{p.v.}
\DeclareMathOperator{\sign}{sign}

\theoremstyle{plain}
\newtheorem{teo}{Theorem}[section]
\newtheorem{newteo}{Theorem}
\newtheorem{cor}[teo]{Corollary}
\newtheorem{lem}[teo]{Lemma}
\newtheorem{proposition}[teo]{Proposition}
\theoremstyle{definition}
\newtheorem{definition}[teo]{Definition}
\newtheorem{oss}[teo]{Remark}

\newtheorem{ex}[teo]{Example}
\newtheorem{ese}[teo]{Exercise}

\numberwithin{equation}{section}

\newcommand{\bpro}{\begin{proposition}}		% Ambienti latex
\newcommand{\npro}{\end{proposition}}		%	|
\newcommand{\bdef}{\begin{definition}}		%	|
\newcommand{\ndef}{\end{definition}}		%	|
\newcommand{\bteo}{\begin{teo}}			%	|
\newcommand{\nteo}{\end{teo}}			%	|
\newcommand{\boss}{\begin{oss}}			%	|
\newcommand{\noss}{\end{oss}}			%	|
\newcommand{\bcor}{\begin{cor}}			%	|
\newcommand{\ncor}{\end{cor}}			%	|
\newcommand{\bdim}{\begin{proof}}			%	|
\newcommand{\ndim}{\end{proof}}			%	|
\newcommand{\bequ}{\begin{equation}}		%	|	
\newcommand{\nequ}{\end{equation}}			%	|
\newcommand{\blem}{\begin{lem}}			%	|
\newcommand{\nlem}{\end{lem}}			%	|
\newcommand{\bex}{\begin{ex}}			%	|
\newcommand{\nex}{\end{ex}}				%	|
\newcommand{\bese}{\begin{ese}}			%	|
\newcommand{\nese}{\end{ese}}			%	|
\newcommand{\N}{\mathbb{N}}				% Spazi matematici
\newcommand{\R}{\mathbb{R}}				%	|
\newcommand{\E}{\mathscr{E}}				%	|
\newcommand{\G}{\mathscr{G}}				%	|
\renewcommand{\L}{\mathcal{L}}			   	%	|
\newcommand{\z}{\textbf{\z}}				%	|
\newcommand{\loc}{{\rm loc}}
\newcommand{\ext}{{\rm ext}}
\title[One-dimensional solutions of non-local Allen-Cahn-type equations]{One-dimensional solutions of non-local Allen-Cahn-type\\ equations with rough kernels}

\author[Matteo Cozzi, Tommaso Passalacqua]{
Matteo Cozzi${}^{(1,2)}$
\and
Tommaso Passalacqua${}^{(1)}$
}

\date{\today}



\begin{document}

\maketitle

{\scriptsize \begin{center} (1) -- Dipartimento di Matematica
``Federigo Enriques''\\
Universit\`a degli Studi di Milano,\\
Via Saldini 50, I-20133 Milano (Italy).\\
\end{center}
\scriptsize \begin{center} (2) -- Laboratoire
Ami\'enois de Math\'ematique Fondamentale et Appliqu\'ee\\
UMR CNRS 7352, Universit\'e de Picardie ``Jules Verne'' \\
33 Rue St Leu, F-80039 Amiens (France).\\
\end{center}
\medskip
\begin{center}
E-mail addresses: matteo.cozzi@unimi.it,
tommaso.passalacqua@unimi.it
\end{center}
}

\bigskip

% Abstract -------------------------------------------------------------------------------------------------------------------------------------------------
	
\begin{abstract}
	We are interested in the study of local and global minimizers for an energy functional of the type
	$$
		\frac{1}{4} \iint_{\R^{2 N} \setminus \left( \R^N \setminus \Omega \right)^2} |u(x) - u(y)|^2 K(x - y) \, dx dy + \int_{\Omega} W(u(x)) \, dx,
	$$
	where~$W$ is a smooth, even double-well potential and~$K$ is a non-negative symmetric kernel in a general class, which contains as a particular case the choice~$K(z) = |z|^{- N - 2 s}$, with~$s \in (0, 1)$, related to the fractional Laplacian. We show the existence and uniqueness (up to translations) of one-dimensional minimizers in the full space~$\R^N$ and obtain sharp estimates for some quantities associated to it. In particular, we deduce the existence of solutions of the non-local Allen-Cahn equation
	$$
		\pv \int_{\R^N} \left( u(x) - u(y) \right) K(x - y) \, dy + W'(u(x)) = 0 \quad \mbox{for any } x \in \R^N,
	$$
	which possess one-dimensional symmetry.
	
	\noindent
	The results presented here were proved in~\cite{CS-M05,PSV13,CS15} for the model case~$K(z) = |z|^{- N - 2 s}$. In our work, we consider instead general kernels which may be possibly non-homogeneous and truncated at infinity.
\end{abstract}

\bigskip

\tableofcontents

%Introduction---------------------------------------------------------------------------------------------------------------------------------------------

\section{Introduction and description of the model}

\noindent
In the present paper we are concerned with a minimization problem related to phase transition phenomena. We study monotone entire minimal configurations for a total energy functional obtained by coupling a standard Gibbs-type free energy with a non-local penalization term modelled upon a Gagliardo-type seminorm. The novelty of our work mostly resides in the introduction of this last term, thanks to which we are able to encompass the presence of long-range interactions between the particles constituting the medium. In particular, our model is general enough to allow for anisotropic effects (possibly changing at different scales of distances, too) and both finite- and infinite-range interactions.

We now proceed to the formal description of the setting.

\medskip

Given a domain~$\Omega \subseteq \R^N$, for some integer~$N \ge 1$, we consider the energy functional
\begin{equation} \label{energy}
	\E_K(u,\Omega):= \mathscr{K}_K(u,\Omega) + \mathscr{P}(u,\Omega),
\end{equation}
where the non-local interaction term~$\mathscr{K}_K$ and the potential term~$\mathscr{P}$ are respectively defined as
$$
	\mathscr{K}_K(u,\Omega) := \frac{1}{4} \iint_{\R^{2N} \setminus (\R^N \setminus \Omega)^2} |u(x) - u(y)|^2 K(x - y) \, dx dy,
$$
and
$$
	\mathscr{P}(u,\Omega):= \int_\Omega W(u(x)) \, dx.
$$

Here,~$K: \R^N \to [0, +\infty]$ is a measurable function modelled on the kernel of the fractional Laplacian. In particular, we ask the kernel~$K$ to fulfill the symmetry condition
\begin{equation} \label{symmetry} \tag{K1}
	K(z) = K(-z) \quad \mbox{for a.a. } z \in \R^N,
\end{equation}
along with various growth and ellipticity assumptions. We highlight the fact that no regularity is required on~$K$. Also notice that, in view of~\eqref{symmetry}, the term~$\mathscr{K}_K$ may be equivalently expressed in the form
$$
	\mathscr{K}_K(u,\Omega) = \frac{1}{4} \int_\Omega \int_{\Omega} |u(x) - u(y)|^2 K(x - y) \, dx dy + \frac{1}{2} \int_\Omega \int_{\R^N \setminus \Omega}|u(x) -u(y)|^2 K(x - y) \,dx dy.
$$

The main ellipticity hypothesis on~$K$ will be\footnote{As it is customary, we denote with~$B_r(x_0)$ the open~$N$-dimensional ball of radius~$r > 0$, centered at a point~$x_0 \in \R^N$. We drop the reference to the center~$x_0$ when~$x_0$ is the origin. That is,~$B_r := B_r(0)$. Moreover, we sometimes write~$S^{N - 1}$ for the unit sphere of~$\R^N$, i.e.~$S^{N - 1} := \partial B_1$. Also,~$\omega_N$ indicates the Lebesgue measure of the~$N$-dimensional ball~$B_1$.}
\begin{equation} \label{weakellipticity} \tag{K2}
	\exists \, 0 < \lambda \le \Lambda, \, r_0 > 0 \, : \, \frac{\lambda \chi_{B_{r_0}(z)}}{|z|^{N+2s}} \le K(z) \le \frac{\Lambda}{|z|^{N+2s}} \quad \mbox{for a.a. } z \in \R^N,
\end{equation}
for some~$s \in (0, 1)$. Notice that condition \eqref{weakellipticity} is very general and allows for a great variety of translation invariant kernels only locally comparable to that of the fractional Laplacian, which is given by the choice~$K(z) = |z|^{- N - 2 s}$. For instance, under~\eqref{weakellipticity} we encompass truncated kernels of the form
\begin{equation} \label{trunkerex}
K(z) = \chi_{B_{r_0}(z)} \frac{a(z)}{|z|^{N + 2 s}},
\end{equation}
with~$a$ bounded and positive, which have been considered in \cite{KKL14,HR-OSV15}. Kernels satisfying~\eqref{weakellipticity}, and even broader similar requirements, are by now widely studied. See e.g.~\cite{K09,K11,DCKP14,DCKP15,CV15}.

For some purposes, we will need the kernel~$K$ to satisfy the stronger condition
\begin{equation} \label{ellipticity} \tag{K2\ensuremath{'}}
	\exists \, 0 < \lambda \le \Lambda \, : \, \frac{\lambda}{|z|^{N+2s}} \leq K(z) \leq \frac{\Lambda}{|z|^{N+2s}} \quad \mbox{for a.a. } z \in \R^N.
\end{equation}
Assumption~\eqref{ellipticity} differs from~\eqref{weakellipticity} in that~$K$ is here required to control the kernel of the fractional Laplacian at all scales and not only in a neighbourhood of the origin. Such hypothesis is more frequently adopted in the literature. To name a few, see~\cite{CS09,CS11,S14,KMS15,KMS15b}.

Finally, to obtain some additional specific results we will restrict ourselves to homogeneous kernels. That is, we will ask~$K$ to be in the form
\begin{equation} \label{Khom} \tag{K2\ensuremath{''}}
K(z) = \frac{a(z / |z|)}{|z|^{N + 2 s}} \quad \mbox{for a.a. } z \in \R^N \mbox{ and with } 0 < \lambda \le a(\zeta) \le \Lambda \mbox{ for a.a. } \zeta \in S^{N - 1}.
\end{equation}
Note that, in dimension~$N = 1$, this and the symmetry condition~\eqref{symmetry} force~$K$ to be the kernel of the fractional Laplacian, up to a multiplicative constant, i.e. 
\begin{equation} \label{N=1Khom}
K(z) = \lambda_\star |z|^{- 1 - 2 s} \quad \mbox{for a.a. } z \in \R,
\end{equation}
for some~$\lambda_\star \in [\lambda, \Lambda]$. We remark that this condition and other generalizations in the same spirit are also often considered in the literature. The interest in~\eqref{Khom} is motivated, for example, by its relationship with stable L\'evy processes in probability theory. On the analysis side, they often lead to slightly sharper results, especially in regularity theory. We refer the interested reader to the works~\cite{FV14,R-OS14b,R-OS15,R-OV15}.

\smallskip

On the other hand, the term~$\mathscr{P}$ is driven by a smooth, even double-well potential~$W$ with wells at~$\pm 1$. More precisely,~$W:\R \to [0, +\infty)$ is a function of class~$C^{2, \beta}_\loc(\R)$, for some~$\beta > 0$, such that
\begin{subequations}
	\begin{align}
		\label{Wpos}  \tag{W1} & W(r) >0 \quad \forall \, r \in (-1,1), \\
		\label{W'=0}  \tag{W2} & W( \pm 1) = W'(\pm 1 ) =0, \\
		\label{W''>0} \tag{W3} & W''(\pm 1 ) >0,
	\end{align}
\end{subequations}
and
\begin{equation} \label{Weven} \tag{W4}
W(r) = W(-r) \quad \forall \, r \in [-1, 1].
\end{equation}
A typical example for~$W$ is represented by the choice
$$
W(r) = \frac{(1 - r^2)^2}{4}.
$$

\smallskip

In this paper we focus on the study of the minimizers for the non-local energy functional~\eqref{energy}. Note that such minimizers are particular solutions of the Euler-Lagrange equation associated to~\eqref{energy}, which is given by
\begin{equation} \label{ELequation}
	-L_Ku + W'(u) = 0,
\end{equation}
where~$L_K$ is the integral operator formally defined as
\begin{equation} \label{operator L form 1}
	-L_K u(x) : = \pv \int_{\R^N} (u(x) - u(y)) K(x - y) \, dy = \lim_{\varepsilon \rightarrow 0^+} \int_{\R^N \setminus B_{\varepsilon}(x)} (u(x) - u(y)) K(x - y) \, dy.
\end{equation}
By changing variables appropriately, we see that~\eqref{operator L form 1} may be equivalently written as
\begin{equation} \label{operator L form 2}
	- L_K u(x) = \pv \int_{\R^N} (u(x) - u(x + z)) K(z) \, dz = \lim_{\varepsilon \rightarrow 0^+} \int_{\R^N \setminus B_{\varepsilon}} (u(x) - u(x + z)) K(z) \, dz.
\end{equation}
Since~\eqref{symmetry} is in force,~$L_K$ can be also represented as a non-singular integral. Indeed, it holds
\begin{equation} \label{operator L form 3}
	L_K u(x) = \frac{1}{2} \int_{\R^N} \delta u(x, z) K(z) \, dz,
\end{equation}
where~$\delta u(x, z)$ is the double increment
\begin{equation} \label{doubleincquot}
	\delta u(x, z) := u(x + z) + u(x - z) - 2u(x).
\end{equation}
We stress that the minus sign in the preceding definitions is chosen so that~$- L_K$ is a positive operator. With this notation, in the special case~$K(z) = |z|^{-N-2s}$ we have that~$- L_K$ is the~$s$-th power of the minus Laplacian, that is
$$
	- L_K u (x) = (-\Delta)^s u(x) = \pv \int_{\R^N} \frac{u(x) - u(y)}{|x -y|^{N+2s}} \, dy,
$$
up to a multiplicative constant. In such situation,~\eqref{ELequation} becomes
\begin{equation} \label{ACfraclap}
(-\Delta)^s u + W'(u) = 0,
\end{equation}
which is often credited as a non-local analogue of the so-called (elliptic)~\emph{Allen-Cahn equation} - the classical, local one being just~\eqref{ACfraclap} with~$s = 1$, formally.

\smallskip

The study of the solutions of the Allen-Cahn equation has been a deep field of research in the last three decades, both in the local and non-local case. Indeed, since the Ginzburg-Landau functional can be viewed as a prototype for the modelling of phase transition phenomena within the Van der Walls-Cahn-Hilliard theory, solutions of the elliptic Allen-Cahn equation represent stationary configurations in this theory.

In the local case, it is well known by the pioneering works of L.~Modica and S.~Mortola~(\cite{MM77}) and E.~De~Giorgi~(\cite{DG79}) that a deep connection between the minimizers of Ginzburg-Landau functionals and minimal surfaces exists. It is probably this relation that prompted De~Giorgi to make his famous conjecture on the symmetry of monotone entire solutions of the Allen-Cahn equation, which eventually paved the way for years of research in nonlinear analysis. See~\cite{BCN97,GG98,AC00,S09,dPKW11} for important contributions in this direction.

\smallskip

In the non-local scenario, there are interesting variations of the above mentioned problems which have attracted the attention of many mathematicians in recent years. An exhaustive report on the various achievements is beyond the scopes of the present work and we instead refer the reader to the surveys~\cite{FV13,BV15}. Nevertheless, we just recall here some of the contributions that are more closely related to the results that will be discussed in the remainder of the paper.

The relationship between the solutions of the fractional Allen-Cahn equation~\eqref{ACfraclap} and minimal surfaces (both the classical ones and an appropriate non-local version of them) is studied in~\cite{SV12}. On the other hand, a suitable fractional version of De~Giorgi conjecture may be stated as follows.
\begin{align*}
& \mbox{Let~$u$ be a bounded entire solution of~\eqref{ACfraclap}, with~$\partial_{x_N} u > 0$ in~$\R^N$.} \\
& \mbox{Is it true that~$u$ must be one-dimensional, i.e. that there exists~$e \in S^{N - 1}$ and~$u_0: \R \to \R$ such that} \\
& \mbox{it holds~$u(x) = u_0(e \cdot x)$ for any~$x \in \R^N$, at least when the dimension~$N$ is~\emph{low}? How low?}
\end{align*}
A positive answer to this question has been given in~\cite{SV09,CS11} for~$N = 2$ and in~\cite{CC10,CC14} for~$N = 3$ and~$s \ge 1/2$. We also report the very recent~\cite{HR-OSV15}, where the authors addressed the validity of such statement in the framework of equation~\eqref{ELequation}, for a class of truncated kernels.

A far more basilar issue in the fractional setting is even the existence itself of one-dimensional solutions. In fact, due to the lack of a satisfactory non-local ODE theory, this problem is not trivial at all. In the case of the fractional Laplacian, it has been solved in~\cite{CS-M05}, for~$s = 1/2$, and in the papers~\cite{PSV13,CS14,CS15}, for a general~$s \in (0, 1)$. We also cite~\cite{AB98,AB98b}, where similar results have been obtained for a class of operators driven by rather general integrable kernels.

\smallskip

In the present work we address precisely this existence result - along with some sharp asymptotic and energy estimates - under hypotheses~\eqref{symmetry} and~\eqref{weakellipticity} (or sometimes~\eqref{ellipticity} and~\eqref{Khom}) on the kernel~$K$. To do this, we follow the lines of the arguments developed in~\cite{PSV13} and suitably adjust them in relation to the changes in our framework. Note that we do not adopt the viewpoint of, say,~\cite{CS15}, as this relies on the so-called Caffarelli-Silvestre extension~(\cite{CS07}), while~\cite{PSV13} does not. This powerful tool enables the interpretation of equations driven by the fractional Laplacian as more common local equations in divergence form. Unfortunately, such extension theory is not available for non-local operators~$L_K$ which differ from the fractional Laplacian. In view of the generality allowed by our setting, we therefore need to undertake a more direct and intrinsically non-local approach.

\smallskip

The following section contains the rigorous statements of our main results.

%--------------------------------------------------------------------------------------------------------------------------------------------------

\section{Main results}
\noindent

As a first step towards the statement of our main contributions, we first need to be precise about the notions of minimizers that we take into consideration. We begin by specifying the definition of local minimizers in bounded domains.

\bdef \label{locmindef}
Let~$\Omega$ be a bounded domain of~$\R^N$. A measurable function~$u: \R^N \to \R$ is said to be a~\emph{local minimizer} for~$\E_K$ in~$\Omega$ if~$\E_K(u, \Omega) < +\infty$ and
$$
\E_K(u, \Omega) \le \E_K(v, \Omega) \quad \mbox{for any measurable } v: \R^N \to \R \mbox{ such that } v = u \mbox{ a.e. in } \R^N \setminus \Omega.
$$
\ndef

As the next remark points out, this concept of local minimization is consistent with respect to set inclusion.

\boss \label{localminsubsetrk}
We observe that if~$\Omega' \subset \Omega$ are two given domains of~$\R^N$, then a local minimizer~$u$ for~$\E_K$ in~$\Omega$ is also a local minimizer in~$\Omega'$. This essentially follows from the facts that
\begin{equation} \label{crossinclusion}
\R^{2 N} \setminus \left( \R^N \setminus \Omega' \right)^2 \subset \R^{2 N} \setminus \left( \R^N \setminus \Omega \right)^2,
\end{equation}
and
$$
\left| u(x) - u(y) \right|^2 = \left| v(x) - v(y) \right|^2 \quad \mbox{for any } (x, y) \in \left( \R^N \setminus \Omega' \right)^2,
$$
if~$u$ and~$v$ coincide outside~$\Omega'$. See~\cite[Remark~1.2]{CV15} for a more detailed explanation of this feature.

Notice that~\eqref{crossinclusion} also implies that the energy~$\E_K(u, \cdot)$ is non-decreasing with respect to set inclusion. In particular, the map
$$
R \longmapsto \E_K(u, B_R),
$$
is monotone non-decreasing, for~$R > 0$.
\noss

We are now in position to provide a satisfactory definition of what a minimizer on the entire space~$\R^N$ is. Note that we can not simply require Definition~\ref{locmindef} to hold with~$\Omega = \R^N$, as the energy~$\E_K$ extended to the full space~$\R^N$ typically diverges. Thus, we shift to the concept of class~A minimizers.

\bdef \label{classAmindef}
A measurable function~$v: \R^N \to \R$ is said to be a~\emph{class~A minimizer} for~$\E_K$ if it is a local minimizer for~$\E_K$ in any bounded domain~$\Omega$ of~$\R^N$.
\ndef

As we just saw, a class~A minimizer is basically a measurable function that minimizes~$\E_K$ with respect to compact perturbations. The terminology we adopted is indeed very classical and tracks back to e.g.~\cite{M24,CdlL01,V04} and, in more recent non-local frameworks close to ours,~\cite{SV14,CV15}.

\smallskip

Our first contribution focuses on the construction of class~A minimizers for~$\E_K$ in one dimension. More precisely we prove the existence and essential uniqueness of a monotone class~A minimizer in the class
\begin{equation} \label{ics}
	\mathcal{X} := \left\{f \in L^1_\loc(\R) \, : \, \lim_{x \to \pm \infty}f(x) = \pm 1\right\},
\end{equation}
of admissible functions. Furthermore, we establish some sharp estimates for the behaviour of such minimizer at infinity and the growth of its energy~$\E_K$ when evaluated on large intervals. To do this, we introduce the quantities
\begin{equation} \label{G*def}
\G_*(u) := \liminf_{R \rightarrow +\infty} \frac{\E_K(u, [-R, R])}{\Psi_s(R)}, \qquad \G^*(u) := \limsup_{R \rightarrow +\infty} \frac{\E_K(u, [-R, R])}{\Psi_s(R)},
\end{equation}
and
$$
\G(u) := \lim_{R \rightarrow +\infty} \frac{\E_K(u, [-R, R])}{\Psi_s(R)},
$$
provided this last limit exists, where
\begin{equation} \label{Psidef}
\Psi_s(R) := \begin{cases}
R^{1 - 2 s} & \quad \mbox{if } s \in (0, 1/2) \\
\log R      & \quad \mbox{if } s = 1/2 \\
1           & \quad \mbox{if } s \in (1/2, 1).
\end{cases}	
\end{equation}
The term~$\Psi_s$ is an important scaling factor that aims at compensating the possible blow up of the energy~$\E_K$ at infinity, in dependence of the parameter~$s$.

The precise statement is as follows.

\begin{newteo} \label{theorem1}
	Let~$N = 1$ and~$s \in (0, 1)$. Assume that~$K$ and~$W$ respectively satisfy conditions~\eqref{symmetry},~\eqref{weakellipticity} and~\eqref{Wpos},~\eqref{W'=0},~\eqref{W''>0},~\eqref{Weven}. Then, there exists an odd, strictly increasing class~A minimizer~$u_0 \in \mathcal{X}$ for~$\E_K$. The minimizer~$u_0$ is of class~$C^{1 + 2 s + \alpha}(\R)$, for some~$\alpha > 0$,\footnote{Note that, given a non-integer~$\gamma > 0$ and a set~$\Omega \subseteq \R^N$, we indicate with~$C^\gamma(\Omega)$ the space composed by functions of~$C^{\lfloor \gamma \rfloor}(\Omega)$ whose partial derivatives of order~$\lfloor \gamma \rfloor$ are globally H\"{o}lder continuous in~$\Omega$, with exponent~$\gamma - \lfloor \gamma \rfloor$. Although no ambiguity should derive from this choice, we will always prefer the more common notation~$C^{\lfloor \gamma \rfloor, \gamma - \lfloor \gamma \rfloor}$ whenever the value of~$\lfloor \gamma \rfloor$ is known.} and is the unique (up to translations) non-decreasing solution\footnote{We point out that by \emph{solution} of an equation like~\eqref{ELequation2} we always mean~\emph{pointwise solution}. We refer the reader to Subsection~\ref{regdefsubsec}, which contains the exact definitions of the two notions of solutions - pointwise and weak - that will be taken into consideration in the paper.} of the Euler-Lagrange equation
	\begin{equation} \label{ELequation2}
	L_K u = W'(u) \quad \mbox{in } \R,
	\end{equation}
	in the class~$\mathcal{X}$.
	
Moreover, there exists a constant~$C \ge 1$ such that the following estimates hold:
	\begin{equation} \label{stima1}
	|u_0(x) - \sign(x)| \le \frac{C}{|x|^{2s}} \quad \mbox{and} \quad |u_0'(x)| \le \frac{C}{|x|^{1 + 2s}} \quad \mbox{for any large } |x|,
	\end{equation}
	\begin{equation}
	\label{N=1tailest}
	\int_{-R}^R \int_{\R \setminus [- R, R]} \left| u_0(x) - u_0(y) \right|^2 K(x - y) \, dx dy \le C R^{1 - 2 s} \quad \mbox{for any large } R > 0,
	\end{equation}
	and
	\begin{equation} \label{N=1enest}
	\G^*(u_0) < +\infty.
	\end{equation}
%	Also, for~$s \in (1/2, 1)$,
%	\begin{equation} \label{s>1/2lim}
%	\mbox{there exists } \, \G(u_0) = \E(u_0, \R) \in (0, +\infty).
%	\end{equation}
	If in addition~$K$ satisfies~\eqref{ellipticity}, then we also have
	\begin{equation} \label{N=1tailestbelow}
	\int_{-R}^R \int_{\R \setminus [- R, R]} \left| u_0(x) - u_0(y) \right|^2 K(x - y) \, dx dy \ge \frac{1}{C} R^{1 - 2 s} \quad \mbox{for any large } R > 0,
	\end{equation}
	and
	\begin{equation} \label{N=1enestbelow}
	\G_*(u_0) > 0.
	\end{equation}
	Finally, if~$s = 1/2$ and~$K$ satisfies~\eqref{Khom} - in its form~\eqref{N=1Khom} -, then
	\begin{equation} \label{barElimit}
 		\G(u_0) = \lim_{R \rightarrow +\infty} \frac{\E_K(u_0, [-R, R])}{\log R} \quad \mbox{exists and is finite},
	\end{equation}
	and it holds
	\begin{equation} \label{barEs=1/2}
	\G(u_0) = \frac{\lambda_\star}{2} \left( \lim_{x \rightarrow +\infty} u_0(x) - \lim_{x \rightarrow -\infty} u_0(x) \right)^2 = 2 \lambda_\star.
	\end{equation}
\end{newteo}

\boss
Observe that the oddness of~$u_0$ is a consequence of the parity assumption~\eqref{Weven}. We stress that, apart from this, such hypothesis on the potential~$W$ is only used at a technical point in Section~\ref{thm1sec}, in order to successfully perform a limiting procedure. We strongly believe that an appropriate adaptation of the arguments contained in~\cite[Sections~3 and~4]{PSV13} may lead to the construction of non-symmetric class~A minimizers, in the absence of~\eqref{Weven}.
\noss

\boss \label{s=1/2limexrk}
Note that, when~\eqref{weakellipticity} is in force with~$s > 1/2$, the existence and finiteness of~$\G(u_0)$ can be easily deduced. Indeed, in such case,
$$
\G(u_0) = \lim_{R \rightarrow +\infty} \E_K(u_0, [-R, R]) = \E_K(u_0, \R),
$$
since the limit exists in view of the monotonicity of the energy (recall Remark~\ref{localminsubsetrk}). Moreover, we also know that~$\G(u_0)$ is finite, thanks to~\eqref{N=1enest}. It is also immediate to check that~$\G(u_0) > 0$, as, otherwise,~$u_0$ would be constant.
\noss

\boss
When~$s = 1/2$, a careful analysis of the proof of~\eqref{barElimit}, provided by Proposition~\ref{Ebars=1/2prop} in Section~\ref{thm1sec}, shows that such conclusion still holds if hypothesis~\eqref{N=1Khom} on~$K$ is replaced by the requirement that
\begin{equation} \label{Katinfty}
\mbox{the limit } K_\infty(z) := \lim_{R \rightarrow +\infty} R^2 K(R z) \mbox{ exists for a.a. } z > 1 \mbox{ and defines a measurable function.}
\end{equation}
We stress that condition~\eqref{Katinfty} is really weaker than~\eqref{N=1Khom}. Indeed,~\eqref{Katinfty} is satisfied for instance by any kernel of the form
$$
K(z) = \frac{\lambda_\star + \sigma(|z|)}{|z|^2},
$$
with~$\lambda_\star > 0$ and~$\sigma: [0, +\infty) \to [0, +\infty)$ measurable, bounded and admitting limit at infinity.

If~$K$ only satisfies~\eqref{Katinfty}, besides~\eqref{symmetry} and~\eqref{weakellipticity}, then~\eqref{barEs=1/2} clearly can not be valid as it is. Nevertheless, by following the proof of Proposition~\ref{Ebars=1/2prop}, it is not hard to see that in such case
\begin{equation} \label{Gu0Kinfty}
\G(u_0) = 2 \int_1^{+\infty} K_\infty(z) \, dz.
\end{equation}
Note that~$K_\infty \in L^1((1, +\infty))$ as a consequence of~\eqref{Katinfty} and~\eqref{ellipticity}. Thus, the right-hand side above is finite.

Furthermore, we point out that~\eqref{barElimit} is trivially satisfied by any truncated kernels, such as for instance those of the form~\eqref{trunkerex}. In this case,~$K_\infty \equiv 0$ and therefore~$\G(u_0) = 0$, in view of~\eqref{Gu0Kinfty}. Besides being interesting on its own, this fact reveals in particular that condition~\eqref{weakellipticity} is not strong enough for~\eqref{N=1enestbelow} to hold, at least for the case~$s = 1/2$.

We believe that an interesting related problem would be to understand whether conclusion~\eqref{barElimit} holds for a larger class of kernels or even for any general~$K$ satisfying~\eqref{symmetry} and~\eqref{weakellipticity}/\eqref{ellipticity}.
\noss

Now that we have established the existence of class~A minimizers on the real line, we can address the problem of how this construction translates to the~$N$-dimensional setting, with~$N \ge 2$. In particular, we shall prove the existence of a~\emph{one-dimensional} class~A minimizer, that is a class~A minimizer for~$\E_K$ in~$\R^N$ that depends only on one single variable, say~$x_N$.

To do this, given a kernel~$K: \R^N \to [0, +\infty]$ satisfying~\eqref{symmetry} and~\eqref{weakellipticity}, let~$k: \R \to [0, +\infty]$ be the kernel defined by\footnote{We reserve the primed notations~$x', y', z'$ for variables in~$\R^{N - 1}$ or, equivalently, in the hyperplane~$\R^{N - 1} \times \{ 0 \}$ of~$\R^N$. Similarly, we often denote with~$B_r'(x_0')$ the open~$(N - 1)$-dimensional ball with radius~$r$ and center~$x_0'$. As for~$N$-dimensional balls,~$B_r'$ stands for the ball centered at the origin.}
\begin{equation} \label{1d kernel}
k(t) := \frac{1}{\varpi} \int_{\R^{N - 1}} K \left( z', \frac{t}{\varpi} \right) \, dz',
\end{equation}
where
\begin{equation} \label{varpidef}
\varpi := \left[ \int_{\R^{N - 1}} \left( 1 + |y'|^2 \right)^{- \frac{N + 2 s}{2}} \, dy' \right]^{- \frac{1}{2 s}}.
\end{equation}
Note that the quantity~$\varpi$ is well-defined and positive (see e.g.~\eqref{finite} for a proof of this fact). The kernel~$k$ is a measurable function which clearly fulfills the symmetry requirement~\eqref{symmetry}, as~$K$ does. Furthermore, it is also easy to see that~$k$ satisfies~\eqref{weakellipticity}. Indeed, by applying the change of variables~$y' := \varpi z' / t$, we compute
\begin{align*}
k(t) & = \frac{|t|^{N - 1}}{\varpi^N} \int_{\R^{N - 1}} K \left( \frac{t}{\varpi} y', \frac{t}{\varpi} \right) \, dy' \ge \frac{|t|^{N - 1}}{\varpi^N} \int_{B'_{\sqrt{\frac{r_0^2 \varpi^2}{t^2} - 1}} } K \left( \frac{t}{\varpi} y', \frac{t}{\varpi} \right) \, dy' \\
& \ge \frac{\lambda \varpi^{2 s}}{|t|^{1 + 2s}} \int_{B'_{\sqrt{\frac{r_0^2 \varpi^2}{t^2} - 1}} } \left( 1 + |y'|^2 \right)^{- \frac{N + 2 s}{2}} \, dy' \ge \frac{\lambda \varpi^{2 s}}{|t|^{1 + 2s}} \int_{B'_1} \left( 1 + |y'|^2 \right)^{- \frac{N + 2 s}{2}} \, dy' = \frac{\tilde{\lambda}}{|t|^{1 + 2 s}},
\end{align*}
for some~$\tilde{\lambda} > 0$, provided~$t < \tilde{r}_0 := \varpi r_0 / \sqrt{2}$. Similarly one checks that the right-hand side inequality in~\eqref{weakellipticity} holds true too. Then, we consider the minimizer~$u_0$ for the energy~$\E_k$ given by Theorem~\ref{theorem1}. We extend it to~$N$-dimensions by setting
\begin{equation} \label{u*def}
u^*(x) := u_0(\varpi x_N) \quad \mbox{for any } x \in \R^N.
\end{equation}

In the next result we show that~$u^*$ is a class~A minimizer for~$\E_K$ and deduce some interesting facts on the asymptotics of the energy~$\E_K(u^*, B_R)$, for~$R > 0$ big.

\begin{newteo} \label{theorem2}
	Let~$N \ge 2$ and~$s \in (0,1)$. Assume that~$K$ and~$W$ respectively satisfy conditions~\eqref{symmetry},~\eqref{weakellipticity} and~\eqref{Wpos},~\eqref{W'=0},~\eqref{W''>0},~\eqref{Weven}. Then, the function~$u^*$ defined in~\eqref{u*def} is a class~A minimizer for~$\E_K$.
	
	Furthermore, the following statements holds true.
	\begin{enumerate}[$\bullet$]
	\item If~$s \in (0, 1/2)$, then there exists a constant~$C \ge 1$ such that
	\begin{equation} \label{extenests<1/2}
	\int_{B_R} \int_{\R^N \setminus B_R} \left| u^*(x) - u^*(y) \right|^2 K(x - y) \, dx dy \le C R^{N - 2 s} \quad \mbox{for any large } R>0.
	\end{equation}
	Also, if~$K$ satisfies~\eqref{ellipticity}, then
	\begin{equation} \label{lowerextenests<1/2}
	\int_{B_R} \int_{\R^N \setminus B_R} \left| u^*(x) - u^*(y) \right|^2 K(x - y) \, dx dy \ge \frac{1}{C} R^{N - 2 s} \quad \mbox{for any large } R>0.
	\end{equation}
	\item If~$s = 1/2$, then
	\begin{equation} \label{totenconvs=1/2}
	\liminf_{R \rightarrow +\infty} \frac{\E_K(u^*, B_R)}{R^{N - 1} \log R} = \frac{\omega_{N - 1}}{\varpi} \, \G_*(u_0), \qquad
	\limsup_{R \rightarrow +\infty} \frac{\E_K(u^*, B_R)}{R^{N - 1} \log R} = \frac{\omega_{N - 1}}{\varpi} \, \G^*(u_0),
	\end{equation}
	and
	\begin{equation} \label{extenvanishs=1/2}
	\lim_{R \rightarrow +\infty} \frac{1}{R^{N - 1} \log R} \int_{B_R} \int_{\R^N \setminus B_R} \left| u^*(x) - u^*(y) \right|^2 K(x - y) \, dx dy = 0.
	\end{equation}
	Also, if~$K$ satisfies~\eqref{Khom}, then the inferior and superior limits in~\eqref{totenconvs=1/2} are equal and it actually holds
	\begin{equation} \label{totenconvs=1/2b}
	\lim_{R \rightarrow +\infty} \frac{\E_K(u^*, B_R)}{R^{N - 1} \log R} = \frac{\omega_{N - 1}}{\varpi} \, \G(u_0) = \frac{\lambda_\star \omega_{N - 1}}{2 \varpi} \left( \lim_{x \rightarrow +\infty} u_0(x) - \lim_{x \rightarrow -\infty} u_0(x) \right)^2 = \frac{2 \lambda_\star \omega_{N - 1}}{\varpi},
	\end{equation}
	with
	\begin{equation} \label{lambdastarthm3}
	\lambda_\star := \varpi^{2 s} \int_{\R^{N - 1}} K \left( y', 1 \right) \, dy'.
	\end{equation}
	\item If~$s > 1/2$, then
	\begin{equation} \label{totenconvs>1/2}
	\lim_{R \rightarrow +\infty} \frac{\E_K(u^*, B_R)}{R^{N - 1}} = \frac{\omega_{N - 1}}{\varpi} \, \E_K(u_0, \R),
	\end{equation}
	and
	\begin{equation} \label{extenvanishs>1/2}
	\lim_{R \rightarrow +\infty} \frac{1}{R^{N - 1}} \int_{B_R} \int_{\R^N \setminus B_R} \left| u^*(x) - u^*(y) \right|^2 K(x - y) \, dx dy = 0.
	\end{equation}
	\end{enumerate}
\end{newteo}

Note that Theorem~\ref{theorem2} is the generalization of~\cite[Theorem~3]{PSV13} to our setting. To prove it, we also extend the techniques of~\cite[Section~5]{PSV13} to rather general integral operators driven by possibly non-homogeneous and truncated kernels (and correct some minor flaws).

The verification of the fact that~$u^*$ is a class~A minimizer is based on the following argument. By Theorem~\ref{theorem1}, we know that the function~$u_0$ is a class~A minimizer for~$\E_k$ and a solution of
$$
L_k u_0 = W'(u_0) \quad \mbox{in } \R.
$$
A simple computation (see~\eqref{Lu*Lu0} in Section~\ref{thm2sec}) then shows that~$u^*$ is a solution of
$$
L_K u^* = W'(u^*) \quad \mbox{in } \R^N.
$$
To obtain that~$u^*$ is actually a class~A minimizer for~$\E_K$, we rely on a general result that connects class~A minimizers and monotone solutions with prescribed limits at infinity in one fixed direction.

\begin{newteo} \label{theorem3}
	Let~$N \ge 1$ and~$s \in (0,1)$. Assume that~$K$ and~$W$ respectively satisfy conditions~\eqref{symmetry},~\eqref{weakellipticity} and~\eqref{W'=0}. Let~$u: \R^N \to (-1, 1)$ be a function of class~$C^{1 + 2 s + \gamma}(\R^N)$, for some~$\gamma > 0$. Suppose that~$u$ is a solution of
	\begin{equation} \label{eullag}
		 L_Ku = W'(u) \quad \mbox{in } \R^N,
	\end{equation}
	which satisfies
	\begin{equation} \label{monotonia}
		\partial_{x_N}u(x) > 0 \quad \forall \, x \in \R^N,
	\end{equation}
	and
	\begin{equation} \label{limitatezza}
		\lim_{x_N \to \pm \infty} u(x',x_N) = \pm 1 \quad \forall \, x' \in \R^{N-1}.
	\end{equation}
	Then,~$u$ is a class~A minimizer for~$\E_K$.
\end{newteo}

We observe that both Theorems~\ref{theorem2} and~\ref{theorem3} are of course still valid if we replace the direction~$e_N$ with a generic direction~$e \in S^{N - 1}$. This can be seen for instance by applying an appropriate rotation in the base space~$\R^N$.

\boss \label{remark1}
	Hypothesis~\eqref{monotonia} may be relaxed to a weak monotonicity assumption. That is, we can replace it with
	\begin{equation} \label{monotonia2}
		\partial_{x_N}u(x) \geq 0 \quad \forall \, x \in \R^N,
	\end{equation}
	without altering the validity of Theorem~\ref{theorem3}. Indeed, it can be shown that if~$u$ satisfies~\eqref{eullag},~\eqref{limitatezza} and~\eqref{monotonia2}, then~$u$ in fact satisfies~\eqref{monotonia}. See Lemma~\ref{rk1pflem} in Subsection~\ref{SCPsubsec} for a proof of this fact.
\noss

\smallskip

The remainder of the paper is organized as follows.

In Section~\ref{regsection} we gather all the regularity results that we will need throughout the exposition. Furthermore, we include there the precise definitions of the functional spaces and notions of solutions that will be used.

In Section~\ref{auxressec} we collect a few preliminary results regarding auxiliary barriers, equations, minimizers and some integral manipulations. In particular, we point the attention of the reader to Subsection~\ref{SCPsubsec}, where we obtain a strong comparison principle for semilinear equations driven by the operator~$L_K$, for a rather general class of non-negative kernels~$K$.

The conclusive three sections are devoted to the proofs of the main results. In Section~\ref{thm3sec} we show the validity of Theorem~\ref{theorem3}. The subsequent Section~\ref{thm1sec} contains the arguments leading to the proof of Theorem~\ref{theorem1}, while the verification of Theorem~\ref{theorem2} occupies the final Section~\ref{thm2sec}.

%--------------------------------------------------------------------------------------------------------------------------------------------------

\section{Regularity of the solutions} \label{regsection}
\noindent
In this section we address the differentiability properties shared by the weak solutions of the linear non-local equation
\begin{equation} \label{Equation}
- L_K u = f \quad \mbox{in } \Omega,
\end{equation}
and of the associated Dirichlet problem
\begin{equation} \label{Dirichlet problem}
	\begin{cases}
      - L_K u = f & \mbox{in } \Omega \\
		    u = g & \mbox{in } \R^N \setminus \Omega,
	\end{cases}
\end{equation}
where~$\Omega$ is a domain of~$\R^N$ and~$f, g$ are measurable functions. Then, we use such results to obtain some informations on the behaviour of the solutions of the semilinear equation~\eqref{ELequation}.

In dependence on how~$\Omega$,~$f$ and~$g$ are chosen, a solution~$u$ may exhibit different regularity features. We do not aim to present here an exhaustive treatise on the regularity theory for~\eqref{Equation}-\eqref{Dirichlet problem} and we instead refer the interested reader to the various contributions available in the literature on the subject (see e.g.~\cite{S06,S07,CS09,CS11,K09,K11,DK12,R-OS14,R-O15,S14}). In fact, we strictly focus on the statements that will be used in the prosecution of the paper.

Furthermore, we point out that the vast majority of the propositions included here are not original and that we intend the present section as a collection of the known regularity results for~\eqref{Equation}-\eqref{Dirichlet problem}, tailored to our needs.

\subsection{Basic definitions} \label{regdefsubsec}

We begin by specifying the notions of solutions that will be adopted throughout the paper. To do this, we first need to introduce the less known functional spaces involved in our definitions. The kernel~$K$ is supposed here to satisfy the general hypotheses~\eqref{symmetry} and~\eqref{weakellipticity}, when not differently stated.

Given any domain~$\Omega \subseteq \R^N$, we consider the linear space
$$
	\mathbb{H}^K(\Omega) := \Big\{  u: \R^N \to \R \mbox{ measurable} : u|_{\Omega} \in L^2(\Omega) \mbox{ and } [u]_{\mathbb{H}^K(\Omega)} < + \infty \Big\},
$$
where
$$
	[u]_{\mathbb{H}^K (\Omega)}^2 := 2 \, \mathscr{K}_K(u, \Omega) = \frac{1}{2} \iint_{\R^{2 N} \setminus \left( \R^N \setminus \Omega \right)^2} \left| u(x) - u(y) \right|^2 K(y-x) \, dx dy.
$$
%Being~$K$ a positive kernel, the space~$H_K(\Omega)$, as endowed with the norm
%$$
%\| u \|_{H_K(\Omega)} := \| u \|_{L^2(\Omega)} + [u]_{H_K(\Omega)},
%$$
%is a Banach space, and even a Hilbert space, with inner product defined by
%$$
%	(u, v)_{H_K(\Omega)}:= \frac{1}{4} \iint_{\R^{2 N} \setminus \left( \R^N \setminus \Omega \right)^2} (u(x) - u(y))(v(x) - v(y)) K(y-x) \, dx dy.
%$$
We point out that
$$
\| u \|_{\mathbb{H}^K(\Omega)} := \| u \|_{L^2(\Omega)} + [ u ]_{\mathbb{H}^K(\Omega)},
$$
is a norm for the space~$\mathbb{H}^K(\Omega)$, as~$K$ is positive near the origin, by~\eqref{weakellipticity}. Moreover, when~$K$ fulfills the stronger condition~\eqref{ellipticity}, then~$\mathbb{H}^K(\Omega)$ is the same as~$\mathbb{H}^s(\Omega)$ - which is just~$\mathbb{H}^k(\Omega)$, with~$k(z) = |z|^{- N - 2 s}$ - considered in~\cite{SerV14}, with equivalent norms. Note that~$\mathbb{H}^K(\Omega)$ differs from the usual fractional Sobolev space~$H^s(\Omega)$ in that the latter does not make any restrictions on the behaviour of its elements outside of~$\Omega$. It holds in fact~$H^s(\R^N) = \mathbb{H}^s(\R^N) \subseteq \mathbb{H}^K(\R^N) \subset \mathbb{H}^K(\Omega) \subset H^s(\Omega)$. Furthermore, we set
$$
\mathbb{H}^K_0(\Omega) := \Big\{ u \in \mathbb{H}^K(\Omega) : u = 0 \mbox{ a.e. in } \R^N \setminus \Omega \Big\} = \Big\{ u \in \mathbb{H}^K(\R^N) : u = 0 \mbox{ a.e. in } \R^N \setminus \Omega \Big\}.
$$

\boss \label{HK0=Hs0rk}
For a general kernel~$K$ satisfying~\eqref{symmetry} and~\eqref{weakellipticity}, it actually holds
\begin{equation} \label{HK0=Hs0}
\mathbb{H}^K_0(\Omega) = \mathbb{H}^s_0(\Omega),
\end{equation}
with equivalent norms, provided~$\Omega$ is bounded. Here,~$\mathbb{H}^s_0(\Omega)$ clearly denotes the subspace of~$\mathbb{H}^s(\Omega)$ composed by the functions vanishing a.e. outside of~$\Omega$.

Notice that, if~\eqref{ellipticity} is in force, then~\eqref{HK0=Hs0} is straightforward. Although not as obvious, the more general assumption~\eqref{weakellipticity} is still strong enough to imply~\eqref{HK0=Hs0}. Indeed, while~\eqref{weakellipticity} ensures that~$K$ and the kernel of the fractional Laplacian are fully comparable only in a neighbourhood of the origin, both these two kernels are integrable at infinity. This and the fact that the functions in~$\mathbb{H}^K_0(\Omega)$ and~$\mathbb{H}^s_0(\Omega)$ are required to vanish outside of~$\Omega$ (the fact that~$\Omega$ has finite measure is of key importance, here) seem to hint at the validity of~\eqref{HK0=Hs0}. Below is a rigorous justification of this quick insight.

First, observe that, by the right-hand inequality in~\eqref{weakellipticity}, it clearly holds~$\mathbb{H}^s_0(\Omega) \subseteq \mathbb{H}^K_0(\Omega)$, with the appropriate inequality for the respective norms. On the other hand, we claim that
\begin{equation} \label{HK0inHs0claim}
[u]_{\mathbb{H}^s(\Omega)} \le c \| u \|_{\mathbb{H}^K(\Omega)} \quad \mbox{for any } u \in \mathbb{H}^K_0(\Omega),
\end{equation}
for some constant~$c > 0$ depending only on~$N$,~$s$,~$\lambda$,~$r_0$ and~$|\Omega|$. Note that, in view of~\eqref{HK0inHs0claim}, equivalence~\eqref{HK0=Hs0} would then follow. Thus, we only need to check~\eqref{HK0inHs0claim}. By using the left-hand side of~\eqref{weakellipticity}, Young's inequality and the fact that~$u = 0$ a.e. in~$\R^N \setminus \Omega$, we compute
\begin{align*}
[u]_{\mathbb{H}^s(\Omega)}^2 & = \frac{1}{2} \int_{\R^N} \left( \int_{B_{r_0}(x)} \frac{\left| u(x) - u(y) \right|^2}{|x - y|^{N + 2 s}} \, dy \right) dx + \frac{1}{2} \int_{\R^N} \left( \int_{\R^N \setminus B_{r_0}(x)} \frac{\left| u(x) - u(y) \right|^2}{|x - y|^{N + 2 s}} \, dy \right) dx \\
& \le \frac{1}{\lambda} [u]_{\mathbb{H}^K(\Omega)}^2 + 2 \left[ \int_{\Omega} |u(x)|^2 \left( \int_{\R^N \setminus B_{r_0}(x)} \frac{dy}{|x - y|^{N + 2 s}} \right) dx + \int_{\Omega} \left( \int_{\Omega \setminus B_{r_0}(x)} \frac{|u(y)|^2}{|x - y|^{N + 2 s}} \, dy \right) dx \right] \\
& \le \frac{1}{\lambda} [u]_{\mathbb{H}^K(\Omega)}^2 + \frac{2}{r_0^{2 s}} \left( N \omega_n + \frac{|\Omega|}{r_0^N} \right) \| u \|_{L^2(\Omega)}^2,
\end{align*}
which is~\eqref{HK0inHs0claim}.
\noss

As a consequence of Remark~\ref{HK0=Hs0rk}, we have that the map
$$
\mathbb{H}^K_0(\Omega) \times \mathbb{H}^K_0(\Omega) \ni (u, v) \longmapsto \langle u, v \rangle_{L^2(\Omega)} + \langle u, v \rangle_{\mathbb{H}^K(\Omega)},
$$
with
\begin{equation} \label{Hilbinner}
\langle u, v \rangle_{\mathbb{H}^K(\Omega)} := \frac{1}{2} \int_{\R^N} \int_{\R^N} \left( u(x) - u(y) \right) \left( v(x) - v(y) \right) K(x - y) \, dx dy,
\end{equation}
is a Hilbert space inner product for~$\mathbb{H}^K_0(\Omega)$, when~$\Omega$ is bounded (see e.g.~\cite[Lemma~7]{SerV12} or~\cite[Lemma~2.3]{FKV15}). Moreover, if~$\Omega$ also has continuous boundary, then
\begin{equation} \label{Cinftydens}
\mathbb{H}^K_0(\Omega) = \overline{C^\infty_0(\Omega)}^{\, \| \cdot \|_{\mathbb{H}^K(\Omega)}},
\end{equation}
as shown in~\cite{FSV15}. We refer to~\cite{DPV12,SerV12,SerV13,FKV15}, to name a few, for additional informations on the above defined spaces and further generalizations.

Throughout the paper we will almost always consider bounded solutions to~\eqref{Equation}. However, for some purposes it is useful to take into consideration a larger class of functions. To this aim, we introduce the weighted Lebesgue space
$$
L^1_s(\R^N) := \Big\{ u: \R^N \to \R \mbox{ measurable } : \| u \|_{L^1_s(\R^N)} < +\infty \Big\},
$$
where
$$
\| u \|_{L^1_s(\R^N)} := \int_{\R^N} \frac{|u(x)|}{1 + |x|^{N + 2 s}} \, dx.
$$
Notice that~$L^\infty(\R^N) \subset L^1_s(\R^N)$, since the weight~$(1 + |x|^{N + 2 s})^{-1}$ is integrable in the whole of~$\R^N$. On the contrary, the space~$L^1_s(\R^N)$ for instance allows for a greater variety of behaviours at infinity.

\smallskip

With all this in hand, we may now head to the definitions of weak solutions of~\eqref{Equation} and~\eqref{Dirichlet problem}.

Let~$\Omega$ be a bounded, Lipschitz domain of~$\R^N$ and~$f \in L^2(\Omega)$. We say that~$u \in \mathbb{H}^K(\Omega)$ is a~\emph{weak solution of equation~\eqref{Equation}} in~$\Omega$ if
\begin{equation} \label{weakform}
\langle u, \varphi \rangle_{\mathbb{H}^K(\Omega)} = \langle f, \varphi \rangle_{L^2(\Omega)} \quad \mbox{for any } \varphi \in \mathbb{H}^K_0(\Omega).
\end{equation}
First, notice that the left-hand side of~\eqref{weakform} is well-defined and finite, as can be seen by inspecting~\eqref{Hilbinner}. Also, in view of~\eqref{Cinftydens}, definition~\eqref{weakform} may be relaxed by requiring it to hold for any~$\varphi \in C^\infty_0(\Omega)$ only, without altering its meaning.

Moreover, given another function~$g \in \mathbb{H}^K(\Omega)$, we say that~$u \in \mathbb{H}^K(\Omega)$ is a~\emph{weak solution of the Dirichlet problem~\eqref{Dirichlet problem}} if~$u - g \in \mathbb{H}^K_0(\Omega)$ and~$u$ weakly solves~\eqref{Equation}.

When~$\Omega$ is not bounded, we may consider a generalized concept of weak solutions of~\eqref{Equation}. In this case,~$u$ is said to be a weak solution of~\eqref{Equation} in~$\Omega$ if, for any Lipschitz subdomain~$\Omega' \subset \subset \Omega$, the function~$u$ belongs to~$\mathbb{H}^K(\Omega')$ and weakly solves~\eqref{Equation} in~$\Omega'$.

When the functions~$u$,~$f$ and~$g$ have more regularity, we may of course strengthen the notion of solution under consideration. Indeed, when~$u \in L^1_s(\R^N) \cap C^{2 s + \gamma}_\loc(\Omega)$, for some~$\gamma > 0$, and~$f$ is, say, continuous in~$\Omega$, then~$u$ is a~\emph{pointwise solution} or, simply, a~\emph{solution} of~\eqref{Equation} if the equation is satisfied at any point~$x \in \Omega$. Similarly, if also~$g \in C^{2 s + \gamma}(\R^N \setminus \Omega)$, then~$u$ is a solution of~\eqref{Dirichlet problem} in~$\Omega$ if~\eqref{Equation} is satisfied in the pointwise sense in~$\Omega$ and~$u \equiv g$ outside of~$\Omega$.

It is immediate to see that~$L_K u(x)$ is well-defined at any point~$x \in \Omega$, when~$u \in L^1_s(\R^N) \cap C^{2 s + \gamma}_\loc(\Omega)$. Also, it is not hard to check that if~$u$ is a weak solution of~\eqref{Equation} and has such regularity, then the equation is also satisfied in the pointwise sense.

\subsection{Linear equations: positive kernels}

In this subsection we enclose all the results that pertain to the linear setting given by~\eqref{Equation}-\eqref{Dirichlet problem}, under the assumption that~$K$ satisfies~\eqref{symmetry} and~\eqref{ellipticity}. In the next subsection, we will remove this latter requirement, by replacing it with the weaker~\eqref{weakellipticity}.

As a first step, we present an interior a priori estimate for the solutions of equation~\eqref{Equation}.

\bpro[\cite{DK12}] \label{uCalphaprop}
Assume that~$K$ satisfies~\eqref{symmetry} and~\eqref{ellipticity}. Let~$f \in L^\infty(B_1)$ and~$u \in L^1_s(\R^N) \cap C^2_\loc(B_1)$ be a solution of~\eqref{Equation} in~$B_1$. Then,~$u \in C^\alpha(B_{1 / 2})$ for any~$\alpha \in (0, \min \{ 2 s, 1 \})$ and it holds
\begin{equation} \label{uCalphaest}
[u]_{C^\alpha(B_{1 / 2})} \le C \left( \| f \|_{L^\infty(B_1)} + \| u \|_{L^1_s(\R^N)} \right),
\end{equation}
for some constant~$C > 0$ which depends only on~$N$,~$s$,~$\lambda$,~$\Lambda$ and~$\alpha$.
\npro

After this preliminary observation, we plan to establish global estimates for the solutions of the Dirichlet problem~\eqref{Dirichlet problem}. For kernels which fulfill the homogeneity condition~\eqref{Khom}, and, actually, more general homogeneous fully nonlinear operators, the optimal~$C^s(\overline{\Omega})$ regularity has been established in~\cite{R-OS15}. In contrast, when~$K$ only satisfies~\eqref{ellipticity}, there is no hope for such boundary regularity, as discussed again in~\cite[Subsection~2.3]{R-OS15}. In the next results we check that it still holds some~$C^\alpha(\overline{\Omega})$ regularity, for~$\alpha < s$.

In conformity with e.g.~\cite{CS09,CS11,R-OS15}, we denote by~$\L_0 = \L_0(s, \lambda, \Lambda)$ the class of operators~$L = L_K$ of the form \eqref{operator L form 3}, whose kernels are measurable functions~$K: \R^N \to [0, +\infty]$ which satisfy \eqref{symmetry} and \eqref{ellipticity}.
The so-called \emph{extremal Pucci operators} for the class~$\L_0$ are defined by
$$
M^+u(x) = M^+_{\L_0}u(x) := \sup_{L \in \L_0} Lu(x) \quad \mbox{and} \quad M^-u(x) = M^-_{\L_0}u(x) := \inf_{L \in \L_0} Lu(x),
$$
For~$\beta \in (0, 2 s)$ and~$\nu \in S^{n - 1}$, we consider the function
$$
\psi_\nu^\beta(x) := \left( \nu \cdot x \right)_+^\beta,
$$
defined for any~$x \in \R^N$.

\bpro[\cite{R-OS15}] \label{CoverCunder}
In correspondence to any~$\beta \in (0, 2 s)$ there exists two constants~$\overline{C}(\beta)$ and~$\underline{C}(\beta)$, which depend on~$N$,~$s$,~$\lambda$ and~$\Lambda$, besides~$\beta$, such that
\begin{align*}
M^+\psi_\nu^\beta(x) & = \overline{C}(\beta) (\nu \cdot x)^{\beta - 2 s} \mbox{ in } \left\{ \nu \cdot x > 0 \right\}, \\
M^-\psi_\nu^\beta(x) & = \underline{C}(\beta) (\nu \cdot x)^{\beta - 2 s} \mbox{ in } \left\{ \nu \cdot x > 0 \right\},
\end{align*}
for every~$\nu \in S^{n - 1}$.

The constants~$\overline{C}$,~$\underline{C}$, viewed as functions of~$\beta$, are continuous in~$(0, 2 s)$. Moreover, there exists two unique values~$0 < \beta_1 < s < \beta_2 < 2 s$, which also depend on~$N$,~$s$,~$\lambda$ and~$\Lambda$, for which
$$
\overline{C}(\beta_1) = 0 = \underline{C}(\beta_2),
$$
and
\begin{align*}
\sign \overline{C}(\beta) & = \sign \left( \beta - \beta_1 \right), \\
\sign \underline{C}(\beta) & = \sign \left( \beta - \beta_2 \right),
\end{align*}
for any~$\beta \in (0, 2 s)$.
\npro

Notice that Proposition~\ref{CoverCunder} is the merging of Proposition~2.7 and Corollary~2.8 in~\cite{R-OS15}. The fact that here the constants~$\overline{C}$ and~$\underline{C}$ do not depend on the direction~$\nu$ is a consequence of the~\emph{isotropy} of the class~$\L_0$. By this we mean that~$\L_0$ is such that
$$
L_K \in \L_0 \quad \mbox{if and only if} \quad L_{K_O} \in \L_0 \mbox{ for any } O \in SO(N),
$$
where~$K_O(z) := K(O z)$. This implies that the Pucci operators~$M^+$ and~$M^-$ are rotationally invariant.\footnote{As noted in~\cite{CS09} the Pucci operators associated to the class~$\L_0$ take the explicit forms
\begin{align*}
M^+u(x) & = \frac{1}{2} \int_{\R^N} \frac{\Lambda \delta u(x, z)_+ - \lambda \delta u(x, z)_-}{|z|^{n + 2 s}} \, dz, \\
M^-u(x) & = \frac{1}{2} \int_{\R^N} \frac{\lambda \delta u(x, z)_+ - \Lambda \delta u(x, z)_-}{|z|^{n + 2 s}} \, dz,
\end{align*}
with~$\delta u(x, z)$ as in~\eqref{doubleincquot}. From this, it is also clear that~$M^+$ and~$M^-$ are rotationally invariant.}

With the aid of the previous proposition, we are now ready to construct a barrier which will eventually prove the {H\"{o}lder continuity of the solutions of~\eqref{Dirichlet problem} up to the boundary of~$\Omega$.

\blem \label{supersollem}
There exist three values~$C \ge 1$,~$r \in (0, 1)$ and~$\beta \in (0, s)$, depending on~$N$,~$s$,~$\lambda$,~$\Lambda$, and a bounded, radial function~$\varphi \in C^{0, \beta}(\R^N) \cap C^\infty(B_{1 + r} \setminus \overline{B_1})$ such that
\begin{equation} \label{supersolprop}
\begin{cases}
M^+ \varphi \le - 1 & \mbox{in } B_{1 + r} \setminus \overline{B_1} \\
\varphi = 0 & \mbox{in } B_1 \\
\varphi(x) \le C \left( |x| - 1 \right)^\beta & \mbox{for any } x \in \R^N \setminus B_1 \\
\varphi \ge 1 & \mbox{in } \R^N \setminus B_{1 + r}.
\end{cases}
\end{equation}
\nlem
\bdim
Let~$\beta_1 \in (0, s)$ be as given by Proposition~\ref{CoverCunder}. Let~$\beta \in (0, \beta_1)$ and define
$$
\varphi^{(\beta)}(x) := \dist \left( x, B_1 \right)^\beta = (|x| - 1)_+^\beta.
$$
We claim that there exists two constants~$\bar{c} > 0$ and~$\bar{r} \in (0, 1)$, depending on~$N$,~$s$,~$\lambda$,~$\Lambda$ and~$\beta$, such that
\begin{equation} \label{prepest1}
M^+ \varphi^{(\beta)}(x) \le - \bar{c} \left( |x| - 1 \right)^{\beta - 2 s} \mbox{ for any } x \in B_{1 + \bar{r}} \setminus \overline{B_1}.
\end{equation}

In order to verify this assertion, we reason as in the proof of Lemma~3.1 in~\cite{R-OS15}. We take~$L = L_K \in \L_0$ and estimate~$L \varphi^{(\beta)}(x_\rho)$, with~$x_\rho = (0, \ldots, 0, 1 + \rho)$ and~$\rho \in (0, 1)$ sufficiently small. To do this, we consider the function
$$
\psi^\beta(x) := \psi_{e_N}^\beta(x - e_N) = \left( x_N - 1 \right)_+^\beta.
$$
It is easy to check that
$$
\psi^\beta(x) \le \varphi^{(\beta)}(x) \mbox{ for any } x \in \R^N, 
$$
and
\begin{equation} \label{preptech2}
\psi^\beta(0, \ldots, 0, x_N) = \varphi^{(\beta)}(0, \ldots, 0, x_N) \mbox{ for any } x_N \in \R.
\end{equation}

By arguing as in the proof of~\cite[Lemma~3.1]{R-OS15}, we also obtain that
$$
\big( \varphi^{(\beta)} - \psi^\beta \big)(x_\rho + z) \le c_1 \begin{cases}
\rho^{\beta - 1} |z'|^2 & \mbox{if } z \in B_{\rho / 2} \\
|z'|^{2 \beta} & \mbox{if } z \in B_1 \setminus B_{\rho / 2} \\
|z|^\beta & \mbox{if } z \in \R^N \setminus B_1, 
\end{cases}
$$
for some constant~$c_1 > 0$. Using this and~\eqref{preptech2}, we estimate
\begin{align*}
L \big( \varphi^{(\beta)} - \psi^\beta \big)(x_\rho) & = \frac{1}{2} \int_{\R^N} \left[ \big( \varphi^{(\beta)} - \psi^\beta \big)(x_\rho + z) + \big( \varphi^{(\beta)} - \psi^\beta \big)(x_\rho - z) \right] K(z) \, dz \\
& \le c_1 \Lambda \left( \int_{B_{\rho / 2}} \frac{\rho^{\beta - 1} |z'|^2}{|z|^{n + 2 s}} \, dz + \int_{B_1 \setminus B_{\rho / 2}} \frac{|z'|^{2 \beta}}{|z|^{n + 2 s}} \, dz + \int_{\R^N \setminus B_1} \frac{|z|^\beta}{|z|^{n + 2 s}} \, dz \right) \\
& \le \frac{c_2}{3} \left( \rho^{\beta + 1 - 2 s} + \rho^{2 \left( \beta - s \right)} + 1 \right) \\
& \le c_2 \rho^{2 \left( \beta - s \right)},
\end{align*}
for some~$c_2 > 0$, since~$\beta < \beta_1 < s$. Thus, recalling Proposition~\ref{CoverCunder}, we get
\begin{align*}
L \varphi^{(\beta)}(x_\rho) & = L \big( \varphi^{(\beta)} - \psi^\beta \big)(x_\rho) + L\psi^\beta(x_\rho) \le c_2 \rho^{2 \left( \beta - s \right)} + M^+ \psi^\beta(x_\rho) \\
& = c_2 \rho^{2 \left( \beta - s \right)} + \overline{C}(\beta) \rho^{\beta - 2 s} = \left( c_2 \rho^\beta - \left| \overline{C}(\beta) \right| \right) \rho^{\beta - 2 s} \\
& \le - \bar{c} \rho^{\beta - 2 s},
\end{align*}
for some~$\bar{c} > 0$, as~$\overline{C}(\beta) < 0$, being~$\beta < \beta_1$, and choosing~$\rho < \bar{r}$, with~$\bar{r} \in (0, 1)$ small enough. Estimate~\eqref{prepest1} then follows by the independence of~$\bar{c}, \bar{r}$ from~$L \in \L_0$ and the rotational symmetry of~$M^+$ and~$\varphi^{(\beta)}$.

Furthermore, if we set
$$
\widetilde{\varphi}^{(\beta)}(x) := \min \left\{ \varphi^{(\beta)}(x), 1 \right\} = \begin{cases}
\left( |x| - 1 \right)_+^\beta & \mbox{if } x \in B_2 \\
1 & \mbox{if } x \in \R^N \setminus B_2,
\end{cases}
$$
then it is not hard to check that
$$
M^+ \widetilde{\varphi}^{(\beta)}(x) \le M^+ \varphi^{(\beta)}(x) + c_3 \le - \bar{c} \left( |x| - 1 \right)^{\beta - 2 s} + c_3 \mbox{ for any } x \in B_{1 + \bar{r}} \setminus \overline{B_1},
$$
for some~$c_3 > 0$. Consequently, by taking a smaller~$\bar{r} > 0$, if necessary, it follows that
$$
M^+ \widetilde{\varphi}^{(\beta)} \le - 1 \mbox{ in } B_{1 + \bar{r}} \setminus \overline{B_1}.
$$

The properties listed in~\eqref{supersolprop} are then satisfied by~$\varphi := C \widetilde{\varphi}^{(\beta)}$, where~$C \ge 1$ is a constant chosen to have~$\varphi \ge 1$ outside of~$B_{1 + \bar{r}}$.
\ndim

Thanks to the supersolution provided by Lemma~\ref{supersollem}, we have

\bpro \label{Calphabdprop}
Assume that~$K$ satisfies~\eqref{symmetry} and~\eqref{ellipticity}. Let~$\Omega \subset \R^N$ be a bounded~$C^{1, 1}$ domain,~$f \in L^\infty(\Omega)$ and~$g \in C^{2 s + \gamma}(\R^N \setminus \Omega)$, with~$\gamma \in (0, 2 - 2 s)$. If~$u \in \mathbb{H}^K(\Omega)$ is a weak solution of the problem~\eqref{Dirichlet problem}, then~$u \in C^\alpha(\overline{\Omega})$, for some~$\alpha \in (0, s)$ depending only on~$N$,~$s$,~$\lambda$,~$\Lambda$ and~$\gamma$, with
\begin{equation} \label{Calphabdest}
\| u \|_{C^\alpha(\overline{\Omega})} \le C \left( \| f \|_{L^\infty(\Omega)} + \| g \|_{C^{2 s + \gamma}(\R^N \setminus \Omega)} + \| u \|_{L^\infty(\Omega)} \right),
\end{equation}
for some constant~$C > 0$ which depends only on~$N$,~$s$,~$\lambda$,~$\Lambda$,~$\gamma$ and~$\Omega$.
\npro

Observe that the we do not need to require a priori the boundedness of~$u$. Indeed, every weak solution of~\eqref{Dirichlet problem} is bounded and satisfies
$$
\| u \|_{L^\infty(\Omega)} \le C (\mbox{diam} (\Omega))^{2 s} \| f \|_{L^\infty(\Omega)},
$$
with~$C > 0$ depending on~$N$,~$s$ and~$\lambda$ (see e.g.~\cite[Corollary~5.2]{R-O15}).

\bdim[Proof of Proposition~\ref{Calphabdprop}]
When~$g \equiv 0$ the proof of~\eqref{Calphabdest} goes as in~\cite[Section~2]{R-OS14}, exploiting Proposition~\ref{uCalphaprop} and Lemma~\ref{supersollem} in place of Corollary~2.5 and Lemma~2.6 there (see~\cite[Subsection~3.4.1]{C16} for more details). The general case then follows by arguing as in~\cite[Remark~7.1]{R-O15}.
\ndim

Next we report a higher order interior regularity result.

\bpro[\cite{S14,R-O15}] \label{C2salphaintprop}
Assume that~$K$ satisfies~\eqref{symmetry} and~\eqref{ellipticity}. Let~$f \in C^\alpha(B_1)$, for some~$\alpha > 0$ such that~$2 s + \alpha$ is not an integer. Let~$u \in \mathbb{H}^K(B_1) \cap C^\alpha(\R^N)$ be a bounded weak solution of~\eqref{Equation} in~$B_1$. Then,~$u \in C^{2 s + \alpha}(B_{1 / 2})$ and
$$
\| u \|_{C^{2 s + \alpha}(B_{1 / 2})} \le C \left( \| f \|_{C^\alpha(B_1)} + \| u \|_{C^\alpha(\R^N)} \right),
$$
for some constant~$C > 0$ which depends only on~$N$,~$s$,~$\lambda$,~$\Lambda$ and~$\alpha$.
\npro

By combining this last result with Proposition~\ref{Calphabdprop}, we obtain the following

\bcor \label{C2salphaintcor}
Assume that~$K$ satisfies~\eqref{symmetry} and~\eqref{ellipticity}. Let~$\Omega \subset \R^N$ be a bounded~$C^{1, 1}$ domain,~$f \in C^{\beta}(\Omega)$ and~$g \in C^{2 s + \gamma}(\R^N \setminus \Omega)$, with~$\beta \in (0, 1)$ and~$\gamma \in (0, 2 - 2 s)$. If~$u \in \mathbb{H}^K(\Omega)$ is a weak solution of~\eqref{Dirichlet problem}, then~$u \in C^{2 s + \alpha}_\loc(\Omega)$, for some~$\alpha \in (0, s)$ depending only on~$N$,~$s$,~$\lambda$,~$\Lambda$,~$\beta$ and~$\gamma$. Also, for any domain~$\Omega' \subset \subset \Omega$ it holds
$$
\| u \|_{C^{2 s + \alpha}(\Omega')} \le C \left( \| f \|_{C^\beta(\Omega)} + \| g \|_{C^{2 s + \gamma}(\R^N \setminus \Omega)} + \| u \|_{L^\infty(\Omega) } \right),
$$
for some constant~$C > 0$ which depends only on~$N$,~$s$,~$\lambda$,~$\Lambda$,~$\beta$,~$\gamma$,~$\Omega$ and~$\Omega'$.
\ncor
%\bdim
%In view of Proposition~\ref{Calphabdprop},~$u \in C^\alpha(\R^N)$, for some~$\alpha \in (0, s)$, and~\eqref{Calphabdest} holds. To obtain that~$u$ is of class~$C^{2 s + \alpha}$ in the interior of~$\Omega$ we refer to~\cite[Theorem~6.1]{R-O15} (see~\cite{S14} for a proof of this fact).
%\ndim

In the next proposition we address the regularity of solutions in the whole space~$\R^N$.

\bpro \label{entirelinregprop}
Assume that~$K$ satisfies~\eqref{symmetry} and~\eqref{ellipticity}. Let~$u \in L^\infty(\R^N)$ be a weak solution of~\eqref{Equation} in~$\R^N$. Then,
\begin{enumerate}[(i)]
\item if~$f \in L^\infty(\R^N)$, then~$u \in C^\alpha(\R^N)$ for any~$\alpha \in (0, \min \{ 2 s, 1 \})$ and
$$
\| u \|_{C^\alpha(\R^N)} \le C \left( \| f \|_{L^\infty(\R^N)} + \| u \|_{L^\infty(\R^N)} \right),
$$
for some constant~$C > 0$ which depends only on~$N$,~$s$,~$\lambda$,~$\Lambda$ and~$\alpha$;
\item if~$f \in C^\alpha(\R^N)$, for some~$\alpha \in (0, 2)$ such that~$2 s + \alpha \ne 1, 2, 3$, then~$u \in C^{2 s + \alpha}(\R^N)$ and
$$
\| u \|_{C^{2 s + \alpha}(\R^N)} \le C \left( \| f \|_{C^\alpha(\R^N)} + \| u \|_{L^\infty(\R^N)} \right),
$$
for some constant~$C > 0$ which depends only on~$N$,~$s$,~$\lambda$,~$\Lambda$ and~$\alpha$.
\end{enumerate}
\npro
\bdim
Item~(i) is an immediate consequence of Proposition~\ref{uCalphaprop} (up to an approximation argument).

On the other hand, to prove~(ii) we first observe that~$u \in C^\beta(\R^N)$, for any~$\beta \in (0, \min \{ 2 s, 1 \})$, in view of~(i). Consider for the moment the case of~$\alpha \in (0, 1)$. If~$s \in (\alpha / 2, 1)$ we may take~$\beta$ to be larger than~$\alpha$. Consequently, both~$u$ and~$f$ belong to~$C^\alpha(\R^N)$ and we are in position to use Proposition~\ref{C2salphaintprop} and recover the~$C^{2 s + \alpha}$ regularity of~$u$.

The case~$s \in (0, \alpha / 2]$ requires a more delicate argument, inspired by an iterative technique displayed in the proof of~\cite[Lemma~6]{PSV13}. Let~$k \ge 1$ be the only integer for which~$s \in (\alpha / (2 k + 2), \alpha / (2 k))]$. Applying Proposition~\ref{C2salphaintprop} for~$k$ times, we get that~$u \in C^{2 k s + \beta}(\R^N)$ for any~$\beta \in (0, 2 s)$, provided~$2 j s + \beta \ne 1$ for each~$j = 1, \ldots, k$. Notice that we are allowed to use this result, since~$\alpha \ge 2 k s > 2 j s + \beta$ for any admissible~$\beta$ and any~$j = 1, \ldots, k - 1$. But then, we can choose~$\beta$ in such a way that~$2 k s + \beta \ge \alpha$, as~$(2 k + 2) s > \alpha$. Hence,~$u \in C^\alpha(\R^N)$ and a further application of Proposition~\ref{C2salphaintprop} leads to the thesis.

When~$\alpha \in [1, 2)$, we already know from the reasoning just displayed that~$u \in C^{2 s + \beta}(\R^N)$ for any~$\beta \in (0, 1)$. Then again, if~$s \in ((\alpha - 1)/2, 1)$, then~$2 s + \beta > \alpha$, for some~$\beta$ close enough to~$1$ and, consequently, we may use Proposition~\ref{C2salphaintprop} to get that~$u \in C^{2 s + \alpha}(\R^N)$. Conversely, when~$s \in (0, (\alpha - 1)/2]$, we argue as before by splitting~$(0, (\alpha - 1) / 2]$ into non-overlapping subintervals. Eventually, we obtain the thesis in this case too.
%
%Under this hypothesis, Proposition~\ref{C2salphaintprop} yields that~$u \in C^{2 s + \beta}(\R^N)$ for any~$\beta \in (0, 2 s)$, provided~$\beta + 2 s \ne 1$. If~$s \ge 1/4$, we may choose~$\beta$ in a way that~$2 s + \beta \ge \alpha$. Then,~$u \in C^\alpha(\R^N)$ and the thesis follows by taking advantage of Proposition~\ref{C2salphaintprop} once again.
%
%Then we are left with~$s \in (0, 1 / 4)$. In this case, by Proposition~\ref{C2salphaintprop},~$u \in C^{4 s + \beta}(\R^N)$ for any~$\beta \in (0, 2 s)$ such that~$\beta + 2 s$
\ndim

We remark that the requirement~$\alpha < 2$ in Proposition~\ref{entirelinregprop}(ii) is only asked for simplicity of exposition. Indeed, one can obtain the result stated there for any~$\alpha > 0$, in the spirit of Proposition~\ref{C2salphaintprop}. However, this formulation is general enough for our future purposes.

\subsection{Linear equations: general kernels}

Here, we extend some results of the previous subsection to operators driven by kernels~$K$ which only satisfy~\eqref{weakellipticity}, instead of the stronger~\eqref{ellipticity}. To do this, we appropriately modify~$K$ far from the origin in order to obtain a new kernel~$\widetilde{K}$ fulfilling~\eqref{ellipticity}. Then, the results will follow by studying the properties of the operator associated to the difference~$\widetilde{K} - K$.

We define~$K_\ext: \R^N \to [0, +\infty)$ to be a radial function of class~$C^\infty$ satisfying
$$
K_\ext(z) = \begin{cases}
\displaystyle 0 & \quad \mbox{if } z \in B_{\frac{r_0}{2}} \\
\displaystyle \frac{\lambda}{|z|^{N + 2 s}} & \quad \mbox{if } z \in \R^N \setminus B_{r_0}.
\end{cases}
$$
The function~$K_\ext$ is clearly bounded. Also, it is not hard to check that~$D^\alpha K_\ext \in L^1(\R^N)$, for every multi-index~$\alpha \in \left( \N \cup \{ 0 \} \right)^N$. In the notation of~\eqref{operator L form 2}, we set
$$
- L_{K_\ext} u(x) := \int_{\R^N} \left( u(x) - u(x - z) \right) K_\ext(z) \, dz.
$$
Observe that~$L_{K_\ext} u$ is well-defined at a.a.~$x \in \R^N$, provided~$u \in L^\infty(\R^N)$. Furthermore,
$$
- L_{K_\ext} u(x) = \| K_\ext \|_{L^1(\R^N)} u(x) - \left( u * K_\ext  \right)(x),
$$
so that~$L_{K_\ext} u$ essentially inherits the regularity properties of~$u$. In particular,
\begin{equation} \label{LinftKext}
\mbox{if } u \in L^\infty(\R^N), \mbox{ then } - L_{K_\ext} u \in L^\infty(\R^N), \mbox{ with } \left\| - L_{K_\ext} u \right\|_{L^\infty(\R^N)} \le C_1 \| u \|_{L^\infty(\R^N)},
\end{equation}
for some constant~$C_1 > 0$ depending on~$K_\ext$, and, given any open set~$\Omega \subseteq \R^N$ and any~$\alpha > 0$,
\begin{equation} \label{CalphaKext}
\mbox{if } u \in L^\infty(\R^N) \cap C^\alpha(\Omega), \mbox{ then } - L_{K_\ext} u \in C^\alpha(\Omega), \mbox{ with } \left\| - L_{K_\ext} u \right\|_{C^\alpha(\Omega)} \le C_2 \| u \|_{C^\alpha(\Omega)},
\end{equation}
for some~$C_2 > 0$ depending on~$K_\ext$ and~$\alpha$.

Let now~$K$ be a kernel satisfying~\eqref{symmetry} and~\eqref{weakellipticity}. We set~$\widetilde{K}(z) := K(z) + K_\ext(z)$, for a.a.~$z \in \R^N$. Notice that the new kernel~$\widetilde{K}$ satisfies~\eqref{symmetry} and~\eqref{ellipticity}, with~$\lambda + \Lambda$ in place of~$\Lambda$. Also,
\begin{equation} \label{HKtilde}
\mbox{if } u \in L^\infty(\R^N) \cap \mathbb{H}^K(\Omega), \mbox{ then } u \in \mathbb{H}^{\widetilde{K}}(\Omega),
\end{equation}
for any bounded domain~$\Omega$.

By knowing all these facts, we are able to extend Proposition~\ref{Calphabdprop} to the case of general kernels satisfying~\eqref{weakellipticity} and obtain a global~$C^\alpha$ regularity result for bounded solutions of the Dirichlet problem~\eqref{Dirichlet problem}.

\bpro \label{Calphabdprop2}
Assume that~$K$ satisfies~\eqref{symmetry} and~\eqref{weakellipticity}. Let~$\Omega \subset \R^N$ be a bounded~$C^{1, 1}$ domain,~$f \in L^\infty(\Omega)$ and~$g \in C^{2 s + \gamma}(\R^N \setminus \Omega)$, with~$\gamma \in (0, 2 - 2 s)$. If~$u \in L^\infty(\R^N) \cap \mathbb{H}^K(\Omega)$ is a weak solution of the problem~\eqref{Dirichlet problem}, then~$u \in C^\alpha(\overline{\Omega})$, for some~$\alpha \in (0, s)$ depending only on~$N$,~$s$,~$\lambda$,~$\Lambda$ and~$\gamma$, with
$$
\| u \|_{C^\alpha(\overline{\Omega})} \le C \left( \| f \|_{L^\infty(\Omega)} + \| g \|_{C^{2 s + \gamma}(\R^N \setminus \Omega)} + \| u \|_{L^\infty(\Omega)} \right),
$$
for some constant~$C > 0$ which depends only on~$N$,~$s$,~$\lambda$,~$\Lambda$,~$r_0$,~$\gamma$ and~$\Omega$.
\npro
\bdim
By~\eqref{HKtilde}, we have that~$u \in L^\infty(\R^N) \cap \mathbb{H}^{\widetilde{K}}(\Omega)$. Moreover,~$u$ is a weak solution of
$$
\begin{cases}
- L_{\widetilde{K}} u = f - L_{K_\ext} u & \quad \mbox{in } \Omega \\
u = g & \quad \mbox{in } \R^N \setminus \Omega.
\end{cases}
$$
Thanks to~\eqref{LinftKext}, the right-hand side~$f - L_{K_\ext} u$ belongs to~$L^\infty(\Omega)$, and the thesis then follows by an application of Proposition~\ref{Calphabdprop}.
\ndim

Similarly, by using~\eqref{CalphaKext} and Corollary~\ref{C2salphaintcor}, we get

\bpro \label{C2salphaintprop2}
Assume that~$K$ satisfies~\eqref{symmetry} and~\eqref{weakellipticity}. Let~$\Omega \subset \R^N$ be a bounded~$C^{1, 1}$ domain,~$f \in C^{\beta}(\Omega)$ and~$g \in C^{2 s + \gamma}(\R^N \setminus \Omega)$, with~$\beta \in (0, 1)$ and~$\gamma \in (0, 2 - 2 s)$. If~$u \in L^\infty(\R^N) \cap \mathbb{H}^K(\Omega)$ is a weak solution of~\eqref{Dirichlet problem}, then~$u \in C^{2 s + \alpha}_\loc(\Omega)$, for some~$\alpha \in (0, s)$ depending only on~$N$,~$s$,~$\lambda$,~$\Lambda$,~$\beta$ and~$\gamma$. Also, for any domain~$\Omega' \subset \subset \Omega$ it holds
$$
\| u \|_{C^{2 s + \alpha}(\Omega')} \le C \left( \| f \|_{C^\beta(\Omega)} + \| g \|_{C^{2 s + \gamma}(\R^N \setminus \Omega)} + \| u \|_{L^\infty(\Omega)} \right),
$$
for some constant~$C > 0$ which depends only on~$N$,~$s$,~$\lambda$,~$\Lambda$,~$r_0$,~$\beta$,~$\gamma$,~$\Omega$ and~$\Omega'$.
\npro

Finally, we extend Proposition~\ref{entirelinregprop} to obtain the following regularity result for entire solutions of~\eqref{Equation}.

\bpro \label{entirelinregprop2}
Assume that~$K$ satisfies~\eqref{symmetry} and~\eqref{weakellipticity}. Let~$u \in L^\infty(\R^N)$ be a weak solution of~\eqref{Equation} in~$\R^N$. Then,
\begin{enumerate}[(i)]
\item if~$f \in L^\infty(\R^N)$, then~$u \in C^\alpha(\R^N)$ for any~$\alpha \in (0, \min \{ 2 s, 1 \})$ and
$$
\| u \|_{C^\alpha(\R^N)} \le C \left( \| f \|_{L^\infty(\R^N)} + \| u \|_{L^\infty(\R^N)} \right),
$$
for some constant~$C > 0$ which depends only on~$N$,~$s$,~$\lambda$,~$\Lambda$,~$r_0$ and~$\alpha$;
\item if~$f \in C^\alpha(\R^N)$, for some~$\alpha \in (0, 2)$ such that~$2 s + \alpha \ne 1, 2, 3$, then~$u \in C^{2 s + \alpha}(\R^N)$ and
$$
\| u \|_{C^{2 s + \alpha}(\R^N)} \le C \left( \| f \|_{C^\alpha(\R^N)} + \| u \|_{L^\infty(\R^N)} \right),
$$
for some constant~$C > 0$ which depends only on~$N$,~$s$,~$\lambda$,~$\Lambda$,~$r_0$ and~$\alpha$.
\end{enumerate}
\npro

\subsection{Semilinear equations}

This conclusive subsection is devoted to a couple of results concerning semilinear equations. These propositions are specifically the ones that are more closely related to the framework in which the paper is set and will be frequently exploited in the following sections. We stress that~$K$ is asked here to satisfy~\eqref{symmetry} and~\eqref{weakellipticity} only.

First is a result for Dirichlet problems in smooth, bounded domains of~$\R^N$.

\bpro \label{C2salphaintsemiprop}
Assume that~$K$ satisfies~\eqref{symmetry} and~\eqref{weakellipticity}. Let~$\Omega \subset \R^N$ be a bounded~$C^{1, 1}$ domain,~$W \in C^{1, \beta}_\loc(\R)$ and~$g \in C^{2 s + \gamma}(\R^N \setminus \Omega)$, with~$\beta \in (0, 1)$ and~$\gamma \in (0, 2 - 2 s)$. If~$u \in L^\infty(\Omega) \cap \mathbb{H}^K(\Omega)$ is a weak solution of
$$
\begin{cases}
L_K u = W'(u) & \quad \mbox{in } \Omega \\
  u = g     & \quad \mbox{in } \R^N \setminus \Omega,
\end{cases}
$$
then~$u \in C^\alpha(\R^N) \cap C^{2 s + \alpha}_\loc(\Omega)$, for some~$\alpha \in (0, s)$ depending only on~$N$,~$s$,~$\lambda$,~$\Lambda$,~$\beta$ and~$\gamma$.
%Also, for any domain~$\Omega' \subset \subset \Omega$ it holds
%$$
%\| u \|_{C^{2 s + \alpha}(\Omega')} \le C \left( \| W'(u) \|_{C^\beta(\Omega)} + \| g \|_{C^{2 s + \gamma}(\R^N \setminus \Omega)} \right),
%$$
%for some constant~$C > 0$ which depends only on~$N$,~$s$,~$\lambda$,~$\Lambda$,~$\beta$,~$\gamma$,~$\Omega$ and~$\Omega'$.
\npro
\bdim
Being~$W'$ continuous and~$u$ bounded, it is clear that the composition~$W'(u)$ is also bounded. In view of this, we may apply Proposition~\ref{Calphabdprop2} to deduce that~$u \in C^{\alpha'}(\overline{\Omega})$, for some~$\alpha' \in (0, s)$. Accordingly,~$u$ is H\"{o}lder continuous in the whole of~$\R^N$. Furthermore,~$W'(u) \in C^{\beta \alpha'}(\Omega)$ and finally Proposition~\ref{C2salphaintprop2} implies that~$u \in C^{2 s + \alpha}_\loc(\Omega)$, with~$\alpha \in (0, s)$.
\ndim

Next, we address the regularity of bounded solutions to semilinear equations in the full space~$\R^N$.

\bpro \label{entiresemiregprop}
Assume that~$K$ satisfies~\eqref{symmetry} and~\eqref{weakellipticity}. Let~$W \in  C^{2, \beta}_\loc(\R)$, for some~$\beta > 0$, and~$u \in L^\infty(\R^N)$ be a weak solution of
$$
L_K u = W'(u) \quad \mbox{in } \R^N.
$$
Then,~$u \in C^{1 + 2 s + \alpha}(\R^N)$, for some~$\alpha > 0$.
\npro
\bdim
We observe that if we show that
\begin{equation} \label{C1alphaclaim}
u \in C^{1, \alpha}(\R^N) \mbox{ for some } \alpha \in (0, \beta],
\end{equation}
then the proof would be over. Indeed, if~$u$ is this regular, then so is~$W'(u)$ and, hence, Proposition~\ref{entirelinregprop2}(ii) implies that~$u \in C^{1 + 2s + \alpha}(\R^N)$.

Thus, we only have to prove~\eqref{C1alphaclaim}. First, we remark that~$W'(u)$ is bounded. Thence, we can use Proposition~\ref{entirelinregprop2}(i) to deduce that~$u$ is of class~$C^{\alpha'}(\R^N)$ for any~$\alpha' \in (0, \min \{ 2 s,  1 \})$. Now we distinguish between the two cases~$s \ge 1 / 2$ and~$s < 1 / 2$.

When~$s \in (1 / 2, 1)$, we have that~$u \in C^{\alpha'}(\R^N)$ for any~$\alpha' \in (0, 1)$. Consequently,~$W'(u) \in C^{\alpha'}(\R^N)$ and we may exploit Proposition~\ref{entirelinregprop2}(ii) to obtain that~$u \in C^{2 s + \alpha'}(\R^N)$ for any such~$\alpha'$, provided~$2 s + \alpha' \ne 2$. Clearly,~\eqref{C1alphaclaim} follows.
%
%On the other hand, if~$s \in [0, 1/2)$ we only know that~$u \in C^{2 s + \alpha'}(\R^N)$ for every~$\alpha' < 2 s$. When~$s \in (1/4, 1/2)$ this is enough to conclude that~$u$ is Lipschitz. Accordingly, so is the composition~$W'(u)$ and, by Proposition~\ref{entirelinregprop}, it follows that~$u \in C^{1, \alpha}(\R^N)$ for any~$\alpha \in (0, 2 s)$.

The case of~$s \in (0, 1/2]$ is slightly more involved. We deal with it by using an approach analogous to the one that we took in the second part of the proof of Proposition~\ref{entirelinregprop2}. Let~$k \ge 1$ be the only integer for which~$s \in (1/(2 k + 2), 1/(2 k)]$. We already know that~$u \in C^{\alpha'}(\R^N)$ for any~$\alpha' \in (0, 2 s)$. Thus, the composition~$W'(u)$ has the same regularity and we may apply Proposition~\ref{entirelinregprop2}(ii) to recover that~$u \in C^{2 s + \alpha'}(\R^N)$, provided~$2 s + \alpha' \ne 1$. By iterating this last step for~$k$ times, we get that~$u \in C^{2 k s + \alpha'}(\R^N)$ for any~$\alpha' \in (0, 2 s)$ such that~$2 j s + \alpha' \ne 1$, for any~$j = 1, \ldots, k$. But now~$2 k s + 2 s > 1$ and thus~\eqref{C1alphaclaim} follows, as we may take~$\alpha'$ as close to~$2 s$ (from below) as we desire.
\ndim

\section{Auxiliary results} \label{auxressec}

\noindent
In this section we include a few preliminary lemmata that will be employed throughout the remainder of the paper to prove the main theorems.

\subsection{Barriers and applications}

Here we construct a couple of useful auxiliary functions that will be needed later on. All the results stated in this subsection are presented without proofs, as their arguments would be almost identical to the ones established in the literature that they generalize. However, we refer the interested reader to~\cite[Subsection~6.2.1]{C16}, where the proofs of all these results are reported in the exact framework of this paper.

We begin by introducing the following barrier.

\blem \label{barrierlem}
Let~$N \ge 1$,~$s \in (0, 1)$ and assume that~$K$ satisfies~\eqref{symmetry} and~\eqref{weakellipticity}. Given any~$\tau > 0$ there exists a constant~$C \ge 1$, which may depend on~$N$,~$s$,~$\Lambda$ and~$\tau$, such that for any~$R \ge C$ we can construct a symmetric radially non-decreasing function
\begin{equation} \label{wbarrange}
w \in C^{1, 1}\left( \R^N, \left[ -1 + C^{-1} R^{- 2 s}, 1 \right] \right),
\end{equation}
with
\begin{equation} \label{wbar=1}
w = 1 \quad \mbox{in } \R^N \setminus B_R,
\end{equation}
which satisfies
\begin{equation} \label{LKwbar}
\left| L_K w(x) \right| \le \tau \left( 1 + w(x) \right),
\end{equation}
and
\begin{equation} \label{1+wbarest}
\frac{1}{C} \left( R + 1 - |x| \right)^{- 2 s} \le 1 + w(x) \le C \left( R + 1 - |x| \right)^{- 2 s},
\end{equation}
for any~$x \in B_R$.
\nlem

Barriers like the one considered in Lemma~\ref{barrierlem} have been first constructed in~\cite{SV14,PSV13} for the fractional Laplacian and in~\cite{CV15b} for more general non-local operators.

In the next result, we obtain another useful barrier in a one-dimensional setting.

\blem \label{etasuperlem}
Let~$N = 1$,~$s \in (0, 1)$ and assume that~$K$ satisfies~\eqref{symmetry} and~\eqref{weakellipticity}. Let~$\eta \in C^2(\R)$ be a positive function such that
$$
\eta(x) = \frac{1}{|x|^{1 + 2 s}} \quad \mbox{for any } x \in \R \setminus (-1, 1).
$$
Then,
$$
L_K \eta \le \Gamma \eta \quad \mbox{in } \R \setminus (-4, 4),
$$
for some constant~$\Gamma \ge 1$ depending only on~$s$,~$\Lambda$ and~$\| \eta \|_{C^2([-1, 1])}$.
\nlem

With the aid of the previous function, one can prove the following bound from above for the decay at infinity of a subsolution of the linear equation
$$
L_K u = \delta u, \quad \mbox{with } \delta > 0,
$$
set on the real line, away from the origin.
%Note that such estimate will be used in Section~\ref{thm2sec} to inspect the asymptotic behaviour of the minimizers in dimension~$N = 1$.

\blem \label{vsubdecay}
Let~$N = 1$,~$s \in (0, 1)$ and assume that~$K$ satisfies~\eqref{symmetry} and~\eqref{weakellipticity}. Let~$R_0, \delta > 0$ be given constants. Let~$v \in C^{2 s + \gamma}(\R)$, for some~$\gamma > 0$, be a bounded function satisfying
$$
L_K v \ge \delta v \quad \mbox{in } \R \setminus [- R_0, R_0].
$$
Then,
$$
v(x) \le \frac{C}{|x|^{1 + 2 s}} \quad \mbox{for any } x \in \R,
$$
for some constant~$C > 0$ possibly depending on~$s$,~$\lambda$,~$\Lambda$,~$R_0$,~$\delta$ and~$\| v \|_{L^\infty(\R)}$.
\nlem

We stress that lemmata~\ref{etasuperlem} and~\ref{vsubdecay} are simple adaptations of, respectively, Lemma~9 and Corollary~4 in~\cite{PSV13}.

\subsection{A strong comparison principle} \label{SCPsubsec}

This subsection focuses on the derivation of a strong comparison principle for semilinear equations. We will heavily rely on such result throughout both Sections~\ref{thm3sec} and~\ref{thm1sec}.

\bpro \label{SCPprop}
Let~$N \ge 1$,~$s \in (0, 1)$ and assume that~$K$ satisfies~\eqref{symmetry} and~\eqref{weakellipticity}. Let~$f_1, f_2: \R^N \times \R \to \R$ be two continuous functions. Let~$\Omega$ be a domain of~$\R^N$ and~$v, w \in L^\infty(\R^N) \cap C^{2 s + \gamma}(\Omega)$, for some~$\gamma > 0$, be such that
$$
\begin{cases}
L_K v \le f_1(\cdot, v) & \quad \mbox{in } \Omega \\
L_K w \ge f_2(\cdot, w) & \quad \mbox{in } \Omega \\
     v \ge w & \quad \mbox{in } \R^N.
\end{cases}
$$
Suppose furthermore that
\begin{equation} \label{f1gef2}
f_1(x, w(x)) \le f_2(x, w(x)) \quad \mbox{for any } x \in \Omega.
\end{equation}
If there exists a point~$x_0 \in \Omega$ at which~$v(x_0) = w(x_0)$, then~$v \equiv w$ in the whole~$\Omega$.
\npro

In the technical hypothesis~\eqref{f1gef2} the two right-hand sides~$f_1$ and~$f_2$ are required to be appropriately ordered on the range of the subsolution~$w$. The conclusion of the proposition is still true if~\eqref{f1gef2} is asked to hold on the range of~$v$, instead. Of course,~\eqref{f1gef2} is clearly satisfied when~$f_1$ and~$f_2$ are the same function.

\bdim[Proof of Proposition~\ref{SCPprop}]
Let~$\varphi := v - w$ and set
$$
\mathcal{Z}_\varphi := \Big\{ x \in \Omega : \varphi(x) = 0 \Big\}.
$$
By assumption, we know that~$\mathcal{Z}_\varphi$ is non-empty, as~$x_0 \in \mathcal{Z}_\varphi$. Moreover,~$\mathcal{Z}_\varphi$ is closed, thanks to the continuity of~$\varphi$ in~$\Omega$. We now claim that~$\mathcal{Z}_\varphi$ is also open. Indeed, let~$\bar{x} \in \mathcal{Z}_\varphi$. Clearly,~$\varphi \ge 0$ in~$\R^N$,~$\varphi(\bar{x}) = 0$ and
$$
L_K \varphi(\bar{x}) = L_K v(\bar{x}) - L_K w(\bar{x}) \le f_1(\bar{x}, v(\bar{x})) - f_2(\bar{x}, w(\bar{x})) = f_1(\bar{x}, w(\bar{x})) - f_2(\bar{x}, w(\bar{x})) \le 0,
$$
in view of~\eqref{f1gef2}. Accordingly,
$$
0 \ge L_K \varphi(\bar{x}) = \frac{1}{2} \int_{\R^N} \left( \varphi(\bar{x} + z) + \varphi(\bar{x} - z) - 2 \varphi(\bar{x}) \right) K(z) \, dz = \frac{1}{2} \int_{\R^N} \left( \varphi(\bar{x} + z) + \varphi(\bar{x} - z) \right) K(z) \, dz \ge 0.
$$
Since, by condition~\eqref{weakellipticity}, the kernel~$K$ is positive in~$B_{r_0}$, we deduce that~$\varphi = 0$ a.a. in~$B_{r_0}(\bar{x})$. That is,~$\Omega \cap B_{r_0} \subseteq \mathcal{Z}_\varphi$. Hence,~$\mathcal{Z}_\varphi$ is open and, by the connectedness of~$\Omega$, we get that~$\mathcal{Z}_\varphi = \Omega$. This concludes the proof.
\ndim

\boss
By inspecting the proof just displayed, we see that the only hypothesis that we really used on~$K$ to deduce the strong comparison principle is its positivity in a small neighbourhood of the origin. This requirement is of course implied by assumption~\eqref{weakellipticity}. But much more different kernels may also enjoy it, such as for instance integrable ones.
\noss

As a first application of Proposition~\ref{SCPprop}, we can now justify the assertion contained in Remark~\ref{remark1}.

\blem \label{rk1pflem}
Let~$N \ge 1$ and~$s \in (0, 1)$. Assume that~$K$ satisfies~\eqref{symmetry} and~\eqref{weakellipticity}. Let~$u \in C^{1 + 2 s + \gamma}(\R^N)$, for some~$\gamma > 0$, be a solution of~\eqref{eullag} which satisfies~\eqref{limitatezza} and~\eqref{monotonia2}. Then,~$u$ also satisfies~\eqref{monotonia}.
\nlem
\bdim
In view of the regularity of~$u$, we may differentiate~\eqref{eullag} in direction~$e_N$ and find that~$\partial_{x_N} u$ solves the equation
\begin{equation} \label{partueq}
L_K \partial_{x_N} u = W''(u) \partial_{x_N} u \quad \mbox{in } \R^N.
\end{equation}
Suppose now by contradiction that there exists~$x_0 \in \R^N$ at which~$\partial_{x_N} u(x_0) = 0$. If this is the case, then by Proposition~\ref{SCPprop} we deduce that~$\partial_{x_N} u = 0$ in the whole of~$\R^N$, which contradicts hypothesis~\eqref{limitatezza}. Note that we can apply such proposition since the function identically equal to~$0$ is another solution of~\eqref{partueq} and~$\partial_{x_N} u \ge 0$, according to~\eqref{monotonia2}. We therefore conclude that~\eqref{monotonia} holds true.
\ndim

\subsection{Existence and stability results} In this subsection we gather a couple of lemmata concerning the existence of local minimizers for~$\E_K$ in a given domain (recall Definition~\ref{locmindef}) and the stability of semilinear equations like~\eqref{ELequation} under locally uniform limits.

We begin with the existence result.

\blem \label{lem2}
	Let~$N \ge 1$ and~$s \in (0, 1)$. Assume that~$K$ and~$W$ respectively satisfy~\eqref{symmetry},~\eqref{weakellipticity} and~\eqref{W'=0}. Let~$\Omega \subset \R^N$ be a bounded Lipschitz domain and~$w_0 : \R^N \to [-1, 1]$ be a measurable function. Suppose that there exists another measurable function~$w$ which coincides with~$w_0$ in~$\R^N \setminus \Omega$ and such that
$$
	\E_K(w, \Omega) < +\infty.
$$
Then, there exists a local minimizer~$u_*: \R^N \to [-1, 1]$ for~$\E_K$ in~$\Omega$ which coincides with~$w_0$ in~$\R^N \setminus \Omega$.
\nlem

\bdim
	Consider a minimizing sequence~$\{ u_j \}_{j \in \N}$, that is
	$$
		\begin{cases}
			\displaystyle u_j = w_0 \quad \mbox{in } \R^N \setminus \Omega \\
			\displaystyle \E_K(u_j, \Omega) \leq \E_K(w, \Omega) \\
			\displaystyle \lim_{j \to +\infty} \E_K(u_j, \Omega) = \inf \left\{ \E_K(v, \Omega) : \, v = w_0 \mbox{ in } \R^N \setminus \Omega \right\} =: \mu.
		\end{cases}
	$$
	Furthermore, by~\eqref{W'=0} we may assume without loss of generality that
	$$
	\left| u_j \right| \le 1 \quad \mbox{in } \R^N,
	$$
	for any~$j \in \N$.	In view of this and~\eqref{weakellipticity}, we compute
	\begin{align*}
		[u_j]_{H^s(\Omega)}^2 & \le \int_{\Omega} \left( \frac{1}{\lambda} \int_{\Omega \cap B_{r_0}(x)} \left| u_j(x) - u_j(y) \right|^2 K(x - y) \, dy + 4 \| u_j \|_{L^\infty(R^N)}^2 \int_{\Omega \setminus B_{r_0}(x)} \frac{dy}{|x - y|^{N + 2 s}} \right) dx \\
		& \le \frac{4}{\lambda} \, \E_K(u_j, \Omega) + \frac{4 |\Omega|^2}{r_0^{N + 2 s}} \| u_j \|_{L^\infty(\R^N)} \\
		& \le c \left( \E_K(w, \Omega) + 1 \right),
	\end{align*}
	for some constant~$c > 0$ independent of~$j$. Hence,~$\{u_j\}$ is bounded in~$H^s(\Omega)$ and then, using e.g.~\cite[Theorem~7.1]{DPV12}, we deduce that~$\{ u_j \}$ converges, up to a subsequence, to some~$u_*$ in~$L^2(\Omega)$ and, thus, a.e. in~$\Omega$. Fatou's Lemma then yields that
	$$
		\E_K(u_*, \Omega) \le \liminf_{j \to +\infty } \E_K(u_j, \Omega) = \mu.
	$$
	This concludes the proof.
\ndim

Secondly, we have the stability lemma.

\blem \label{stablem}
Let~$N \ge 1$,~$s \in (0, 1)$ and assume that~$K$ satisfies conditions~\eqref{symmetry} and~\eqref{weakellipticity}. Let~$\Omega \subset \R^N$ be a bounded Lipschitz domain and~$\{ v_j \}_{j \in \N} \subset \mathbb{H}^K(\Omega) \cap L^\infty(\R^N)$ a sequence of functions. Assume that~$v_j$ is a weak solution of
\begin{equation} \label{vjstabsol}
L_K v_j = W'(v_j) \quad \mbox{in } \Omega,
\end{equation}
and that there exists a constant~$C > 0$ such that
\begin{equation} \label{vjunif}
[v_j]_{\mathbb{H}^K(\Omega)} + \| v_j \|_{L^\infty(\R^N)} \le C,
\end{equation}
for any~$j \in \N$. Suppose furthermore that~$v_j$ converges to a function~$v$ uniformly on compact subsets of~$\R^N$. Then,~$v \in \mathbb{H}^K(\Omega) \cap L^\infty(\R^N)$ and is a weak solution of
\begin{equation} \label{vstabsol}
L_K v = W'(v) \quad \mbox{in } \Omega.
\end{equation}
\nlem
\bdim
First of all, it is clear that~$v$ belongs to~$L^\infty(\R^N)$, as~$v_j \rightarrow v$ locally uniformly in~$\R^N$ and~\eqref{vjunif} holds. It is immediate to check that~$v \in \mathbb{H}^K(\Omega)$, since, by~\eqref{vjunif} and Fatou's lemma,
\begin{align*}
[v]_{\mathbb{H}^K(\Omega)}^2 & = \frac{1}{2} \iint_{\R^{2 N} \setminus \left( \R^N \setminus \Omega \right)^2} \left| v(x) - v(y) \right|^2 K(x - y) \, dx dy \\
& \le \frac{1}{2} \liminf_{j \rightarrow +\infty} \iint_{\R^{2 N} \setminus \left( \R^N \setminus \Omega \right)^2} \left| v_j(x) - v_j(y) \right|^2 K(x - y) \, dx dy \\
& = \liminf_{j \rightarrow +\infty} \, [v_j]_{\mathbb{H}^K(\Omega)}^2 \\
& \le C^2.
\end{align*}

Now we show that~$v$ is a weak solution of~\eqref{vstabsol}. Fix~$\varphi \in C_0^\infty(\Omega)$. Since~$v_j$ is a weak solution of~\eqref{vjstabsol}, we have that
\begin{equation} \label{vjstabtech}
\begin{aligned}
\int_{\Omega} W'(v_j(x)) \varphi(x) \, dx & = - \frac{1}{2} \int_{\R^N} \int_{\R^N} \left( v_j(x) - v_j(y) \right) \left( \varphi(x) - \varphi(y) \right) K(x - y) \, dx dy \\
& = \int_{\Omega} v_j(x) L_K \varphi(x) \, dx.
\end{aligned}
\end{equation}
Notice now that~$L_K \varphi \in L^\infty(\Omega) \subset L^1(\Omega)$. Indeed, by~\eqref{ellipticity} we have
\begin{align*}
\left| L_K \varphi (x) \right| & = \frac{1}{2} \left| \int_{\R^N} \left[ \varphi(x + z) + \varphi(x - z) - 2 \varphi(x) \right] K(z) \, dz \right| \\
& \le \frac{\Lambda}{2} \left[ 4 \| \varphi \|_{L^\infty}(\R^N) \int_{B_1} \frac{dz}{|z|^{N + 2 s}} + \| \nabla^2 \varphi \|_{L^\infty(\R^N)} \int_{\R^N \setminus B_1} \frac{dz}{|z|^{N - 2 + 2 s}} \right] \\
& \le N \omega_N \Lambda \| \varphi \|_{C^2(\R^N)} \left[ \frac{1}{s} + \frac{1}{1 - s} \right],
\end{align*}
for any~$x \in \Omega$. Hence, by the Dominated Convergence Theorem and the continuity of~$W'$, we may take the limit as~$j \rightarrow +\infty$ in~\eqref{vjstabtech} and deduce that
$$
\int_{\Omega} v(x) L_K \varphi(x) \, dx = \int_{\Omega} W'(v(x)) \varphi(x) \, dx.
$$
Since we have already showed that~$v \in \mathbb{H}^K(\Omega)$, it easily follows that~$v$ is a weak solution of~\eqref{vstabsol}.
\ndim

\subsection{Some integral computations}

We conclude the section with a couple of results aimed at establishing an upper bound for the quantity
\begin{equation} \label{Jdef}
J_{\alpha, N}(\rho, \sigma) := \int_{B_\rho} \int_{\R^N \setminus B_\sigma} {\left( 1 + |x - y|^2 \right)}^{-\alpha} \, dx dy, \qquad \mbox{for } \rho, \sigma > 0 \mbox{ and } \alpha > \frac{N}{2}.
\end{equation}
This will play an important role later in Section~\ref{thm2sec}, to perform some computations needed for the proof of Theorem~\ref{theorem2}.

First, we have the following

\blem \label{Ilem}
Let~$N \ge 1$,~$\alpha \in (N/2,+\infty)$ and~$\rho > \sigma > 0$. Then, given any~$\delta \in (0, 1)$ satisfying
\begin{equation} \label{deltabounds}
\frac{N}{2\alpha} < \delta <  \frac{N + 1}{2\alpha},
\end{equation}
it holds
\begin{equation} \label{Iest1}
\int_{B_\rho} \int_{B_\rho \setminus B_\sigma} (1 + |x - y|^2)^{-\alpha} \, dx dy \le C_1 \left( \rho^N - \sigma^N \right),
\end{equation}
and
\begin{equation} \label{Iest2}
\int_{B_\rho} \int_{\R^N \setminus B_\rho} (1 + |x - y|^2)^{-\alpha} \, dx dy \le C_2 \rho^{2 N - 2 \delta \alpha},
\end{equation}
for some constants~$C_1 > 0$, which depends only on~$N$ and~$\alpha$, and~$C_2 > 0$, which may also depend on~$\delta$.
\nlem

\bdim
All along the proof,~$c$ will denote any positive constant depending on~$N$ and~$\alpha$, whose value may change from line to line.

We begin by establishing~\eqref{Iest1}. Changing variables appropriately we compute
\begin{align*}
\int_{B_\rho} \int_{B_\rho \setminus B_\sigma} (1 + |x - y|^2)^{-\alpha} \, dx dy & \le \int_{B_\rho \setminus B_\sigma} \left( \int_{\R^N} (1 + |z|^2)^{-\alpha} \, dz \right) dx \\
& = |B_\rho \setminus B_\sigma| \int_{\R^N} (1 + |z|^2)^{-\alpha} \, dz \\
& = c_1 \left( \rho^N - \sigma^N \right),
\end{align*}
for some constant~$c_1 > 0$ depending on~$N$ and~$\alpha$. This is true since
\begin{equation} \label{finite}
\begin{aligned}
\int_{\R^N} (1 + |z|^2)^{-\alpha} \, dz & = N \omega_N \int_0^{+\infty} \left( 1 + r^2 \right)^{-\alpha} r^{N - 1} \, dr \\
& \le N \omega_N \left( \int_0^1 r^{N - 1} \, dr + \int_1^{+\infty} r^{N - 2 \alpha - 1} \, dr \right) \\
& = N \omega_N \left( \frac{1}{N} - \frac{1}{N - 2 \alpha} \right) \\
& < +\infty,
\end{aligned}
\end{equation}
as~$N - 2 \alpha < 0$. Therefore,~\eqref{Iest1} is proved.

We now address~\eqref{Iest2}. Consider any real number~$0 < \delta < 1$. From now on,~$c$ is allowed to depend on~$\delta$ too. Applying Young's inequality with weight~$\delta$, we get
$$
1 + |x - y|^2 = (1 - \delta) {\left( \frac{1}{(1 - \delta)^{1 - \delta}} \right)}^{1 / (1 - \delta)} + \delta {\left( \frac{|x - y|^{2 \delta}}{\delta^\delta} \right)}^{1 / \delta} \ge \frac{|x - y|^{2 \delta}}{(1 - \delta)^{1 - \delta} \delta^\delta}.
$$
We estimate
\begin{equation} \label{deltatech}
\begin{aligned}
\int_{B_\rho} \int_{\R^N \setminus B_\rho} (1 + |x - y|^2)^{-\alpha} \, dx dy & \le c \int_{B_\rho} \left( \int_{\R^N \setminus B_{\rho - |y|}(y)} |x - y|^{- 2 \delta \alpha} \, dx \right) dy \\
& = c \int_{B_\rho} \left( \int_{\rho - |y|}^{+\infty} r^{N - 1 - 2 \delta \alpha} \, dr \right) dy.
\end{aligned}
\end{equation}
Now, we require~$\delta$ to satisfy~\eqref{deltabounds}. Under this restriction,~$N - 2 \delta \alpha < 0$ and thus~\eqref{deltatech} becomes
$$
\int_{B_\rho} \int_{\R^N \setminus B_\rho} (1 + |x - y|^2)^{- \alpha} \, dx dy \le c \int_{B_\rho} \left( \rho - |y| \right)^{N - 2 \delta \alpha} \, dy \le c \rho^{N - 1} \int_0^\rho \left( \rho - r \right)^{N - 2 \delta \alpha} \, dr.
$$
But then,~\eqref{deltabounds} also implies that~$N - 2 \delta \alpha + 1 > 0$, so that
$$
\int_{B_\rho} \int_{B_{2 \rho} \setminus B_\rho} (1 + |x - y|^2)^{-\alpha} \, dx dy \le c \rho^{2 N - 2 \delta \alpha},
$$
which is~\eqref{Iest2}.
\ndim

From Lemma~\ref{Ilem} we immediately get the desired estimate for~$J$.

\bcor \label{Jcor}
Let~$N \ge 1$,~$\alpha \in (N/2,+\infty)$ and~$\rho, \sigma > 0$. Then, given any~$\delta \in (0, 1)$ satisfying~\eqref{deltabounds},
$$
J_{\alpha, N}(\rho, \sigma) \le C \Big( \rho^{2 N - 2 \delta \alpha} + \max \left\{ \rho^N - \sigma^N, 0 \right\} \Big),
$$
for some constant~$C > 0$ which depends on~$N$,~$\alpha$ and~$\delta$.
\ncor

% Proof of Theorem 1 -------------------------------------------------------------------------------------------------------------------

\section{Proof of Theorem~\ref{theorem3}} \label{thm3sec}

\noindent
In this section we present a proof of Theorem~\ref{theorem3}. We stress that the argument displayed is an adaptation of that of~\cite[Theorem 1]{PSV13}, in accordance with the changes in our setting.

\smallskip

\emph{Step 1.} Arguing by contradiction, we suppose that~$u$ is not a class~A minimizer for~$\E_K$. Recalling Definition~\ref{classAmindef}, there exists a bounded domain~$\Omega \subset \R^N$ in which~$u$ is not a local minimizer. According to Remark~\ref{localminsubsetrk}, we may further assume that~$\Omega = B_R$, for some~$R >0$. Thus, there exists a function~$\varphi$ supported in~$B_R$ such that
$$
	\E_K(u + \varphi, B_R) < \E_K(u,B_R).
$$
Note that this implies in particular that~$\E_K(u + \varphi, B_R)$ is finite. Hence, we may apply Lemma~\ref{lem2} with~$w = u + \varphi$ and find a minimizer~$u_*$ for~$\E_K(\cdot, B_R)$ among all functions~$v$ such that~$v = u$ outside of~$B_R$. Observe that Lemma~\ref{lem2} also tells us that
$$
|u_*| \le 1 \quad \mbox{in } \R^N.
$$
Since we assumed by contradiction that~$u$ is not a minimizer, there exists a point~$x_0 \in \R^N$ at which~$u_*(x_0) \neq u(x_0)$. We suppose in fact that
\begin{equation} \label{passo1}
	u_*(x_0) > u(x_0).
\end{equation}
A specular argument can be provided in case the opposite inequality holds. By the minimizing property of~$u_*$ we have that~$u_*$ is a weak solution of
\begin{equation} \label{minprop}
	 L_K u_* = W'(u_*) \quad \mbox{in } B_R.
\end{equation}
Therefore, we may apply Proposition~\ref{C2salphaintsemiprop} to conclude that~$u_*$ is continuous in the whole of~$\R^N$. Also, observe that, by the same proposition,~$u_*$ is of class~$C^{2 s + \alpha}$ in the interior of~$B_R$ and thus~\eqref{minprop} holds in the pointwise sense.

\emph{Step 2.} Now we can prove that 
\begin{equation} \label{passo2}
	\left| u_* \right| < 1,
\end{equation}
using the assumptions on the potential~$W$.

Indeed, suppose that there exists~$\bar{x} \in \R^N$ at which, e.g.,~$u_*(\bar{x}) = -1$. Since~$|u| < 1$ and~$u_*$ coincides with~$u$ outside~$B_R$ we conclude that~$\bar{x} \in B_R$. Hence, by also recalling~\eqref{minprop} and~\eqref{W'=0}, we are in position to apply Proposition~\ref{SCPprop} (with~$v = u_*$ and~$w = -1$) to deduce that~$u_* \equiv -1$ in~$B_R$. But this and the continuity of~$u_*$ up to the boundary of~$B_R$ contradict the assumption that~$u_* \equiv u$ outside~$B_R$, as~$|u| < 1$. Then~\eqref{passo2} holds true.

\emph{Step 3.} We claim that there exists~$\bar{k} \in \R$ such that
\begin{equation} \label{passo3}
	\mbox{if } k \geq \bar{k}, \mbox{ then } u(x + k e_N) \geq u_*(x) \quad \forall \, x \in \R^N.
\end{equation}
Again we argue by contradiction and suppose that there exist two sequences~$k_j > 0$ and~$x^{(j)} \in \R^N$ such that~$k_j \rightarrow +\infty$ as~$j \rightarrow +\infty$ and
\begin{equation} \label{u<u*}
	u(x^{(j)} + k_j e_N) < u_*(x^{(j)}).
\end{equation}
Since~$u$ is monotone in the~$e_N$ direction by assumption~\eqref{monotonia} and~$k_j \geq 0$, it follows that
$$
	u(x^{(j)}) < u_*(x^{(j)}),
$$
and therefore~$x^{(j)} \in B_R$. Hence, up to a subsequence,~$x^{(j)}$ converges to some~$x_* \in \overline{B_R}$. But now, taking advantage of assumption~\eqref{limitatezza}, inequality~\eqref{u<u*} and the continuity of~$u_*$ in~$\overline{B_R}$, we find
$$
	1 = \lim_{j \rightarrow +\infty} u(x^{(j)} + k_j e_N) \leq \lim_{j \rightarrow +\infty} u_*(x^{(j)})= u_*(x_*).
$$ 
But this is in contradiction with~\eqref{passo2} and so~\eqref{passo3} is proved.

\emph{Step 4.} Now we can take~$\hat{k}$ as the least possible value of~$\bar{k}$ for which~\eqref{passo3} holds. Thus, there exist two sequences~$\eta_j >0$ and~$y^{(j)} \in \R^N$ for which
\begin{equation} \label{monopt}
	u(y^{(j)} + (\hat{k}- \eta_j) e_N) \leq u_*(y^{(j)}),
\end{equation}
and~$\eta_j \rightarrow 0^+$ as~$j \rightarrow +\infty$. Now, by~\eqref{passo1} and~\eqref{passo3} we have that 
$$
	u(x_0) < u_*(x_0) \leq u(x_0 + \hat{k} e_N),
$$
so that
\begin{equation} \label{hatk>0}
\hat{k} >0,
\end{equation}
by the monotonicity of~$u$.

We claim that there exists~$J \in \N$ such that
\begin{equation} \label{passo4}
	y^{(j)} \in B_R \quad \forall \, j \ge J.
\end{equation}
By contradiction, if~$y^{(j)} \in \R^N \setminus B_R$ for infinitely many~$j$'s, by~\eqref{monopt} and the fact that~$u_* \equiv u$ outside of~$B_R$, we would have that
$$
	u(y^{(j)} + (\hat{k} - \eta_j) e_N) \leq u_*(y^{(j)}) = u(y^{(j)}).
$$
But then 
$$
	\hat{k}- \eta_j \leq 0,
$$
by the monotonicity of~$u$, and thus, by letting~$j$ go to~$+\infty$, we would get~$\hat{k} \leq 0$. But this is contradicts~\eqref{hatk>0} and hence~\eqref{passo4} holds true.

\emph{Step 5.} In view of the previous deduction, we can assume that
$$
	\lim_{j \to +\infty} y^{(j)} = y_*,
$$
for some~$y_*$ in the closure of~$B_R$. Taking the limit as~$j \rightarrow +\infty$ in~\eqref{monopt} and recalling~\eqref{passo3}, we then get
\begin{equation} \label{uy*hatk}
u(y_* + \hat{k} e_N) = u_*(y_*).
\end{equation}
But using once again the strict monotonicity of~$u$ and recalling~\eqref{hatk>0}, we are led to
$$
u(y_*) < u_*(y_*).
$$
Consequently,~$y_* \in B_R$, as~$u$ and~$u_*$ coincide outside of~$B_R$.

Define now~$v(x) := u(x + \hat{k} e_N)$, for any~$x \in \R^N$. By~\eqref{eullag},~\eqref{minprop} and~\eqref{passo3}, we know that
$$
	\begin{cases}
		L_K v = W'(v) & \quad \mbox{in } \R^N \\
	L_K u_* = W'(u_*) & \quad \mbox{in } B_R \\
		    v \ge u_* & \quad \mbox{in } \R^N.
	\end{cases}
$$
Also, by~\eqref{uy*hatk}, we have that~$v(y_*) = u_*(y_*)$. Thence, by applying Proposition~\ref{SCPprop} (with~$w = u_*$ and~$\Omega = B_R$) we obtain that~$v \equiv u_*$ in the whole~$B_R$.

The strict monotonicity of~$u$,~\eqref{hatk>0} and the continuity of~$u$ and~$u_*$ up to the boundary of~$B_R$ imply in turn that
$$
	u(x) < u(x + \hat{k} e_N) = u_*(x) \quad \forall \,  x \in \overline{B_R},
$$
contradicting the fact that~$u$ coincides with~$u_*$ outside~$B_R$. Thus, the proof is complete.

\section{Proof of Theorem~\ref{theorem1}} \label{thm1sec}

\noindent
Here we show the existence of a class~A minimizer in dimension~$N = 1$, thus proving Theorem~\ref{theorem1}. To do so, we first deal with a constraint minimization problem on intervals, in Subsection~\ref{minintsubsec}. Then, in Subsection~\ref{lineminsubsec}, we obtain the existence of local minimizers on the whole real line~$\R$. Finally, the conclusive Subsections~\ref{est1subsec},~\ref{est2subsec} and~\ref{est3subsec} are devoted to the study of the various estimates involved in the statement of Theorem~\ref{theorem1}.

\subsection{Minimizers on intervals} \label{minintsubsec}

In the first proposition of the subsection we deal with the existence of local minimizers on~\emph{large} real intervals and prove some key estimates for their energies.

\blem \label{intminlem}
Let~$M > 3$. Then, there exists a local minimizer~$v_{[-M, M]}: \R \to [-1, 1]$ for~$\E_K$ in~$[-M, M]$, such that~$v_{[-M, M]}(x) = - 1$ for any~$x \le -M$ and~$v_{[-M, M]}(x) = 1$ for any~$x \ge M$. Moreover,~$v_{[-M, M]}$ is odd, non-decreasing and is the unique solution of the Dirichlet problem
\begin{equation} \label{v-MMDir}
\begin{cases}
L_K u = W'(u) & \quad \mbox{in } (-M, M) \\
u = -1        & \quad \mbox{in } (-\infty, -M] \\
u = 1         & \quad \mbox{in } [M, +\infty).
\end{cases}
\end{equation}
Finally, there exists a constant~$C \ge 1$, depending only on~$s$,~$\Lambda$ and~$W$, for which
\begin{equation} \label{vIest}
\E_K(v_{[-M, M]}, J) \le C \Psi_s(|J|),
\end{equation}
where~$J$ is either~$[-M, M]$ or any subinterval of~$[-M, M]$ such that~$|J| > 6$ and~$\dist(J, \R \setminus (-M, M)) > 2$.
\nlem

Recall that the quantity~$\Psi_s$ was defined in~\eqref{Psidef}.

\bdim[Proof of Lemma~\ref{intminlem}]
Consider the piecewise linear function~$h: \R \to [-1, 1]$ defined by
$$
h(x) := \begin{cases}
-1 & \quad \mbox{if } x \le - 1 \\
x  & \quad \mbox{if } -1 < x \le 1 \\
1  & \quad \mbox{if } x > 1.
\end{cases}
$$
By arguing as in~\cite[Lemma~2]{PSV13} and taking advantage of the right-hand inequality in~\eqref{weakellipticity}, it is easy to check that
\begin{equation} \label{hest}
\E_K(h, [-M, M]) \le c \Psi_s(M) < +\infty,
\end{equation}
for some constant~$c > 0$ only depending on~$s$,~$\Lambda$ and~$W$. The existence of a local minimizer~$v_{[-M, M]}: \R \to [-1, 1]$ then derives from Lemma~\ref{lem2}. Note that~\eqref{hest} also establishes~\eqref{vIest} for~$J = [-M, M]$. Estimate~\eqref{vIest} for a general interval~$J \subset [-M, M]]$ with~$|J| > 6$ and~$\dist(J, \R \setminus (-M, M)) > 2$ follows from e.g.~\cite[Proposition~3.1]{CV15},\footnote{Note that Proposition~3.1 in~\cite{CV15} is proved under the analogous of assumption~\eqref{weakellipticity} here, with~$r_0 = 1$. However, the proof of that result only exploits the right-hand inequality of~\eqref{weakellipticity} and thus it is valid in our framework too.} by observing that~$v_{[-M, M]}$ is also a local minimizer in any subinterval of~$[-M, M]$ (recall Remark~\ref{localminsubsetrk}).

Then, we address the monotonicity of~$v_{[-M, M]}$. First, note that~$v_{[-M, M]}$ weakly solves~\eqref{v-MMDir} and, therefore, by Proposition~\ref{C2salphaintsemiprop},~$v_{[-M, M]} \in C^\alpha(\R) \cap C^{2 s + \alpha}_\loc((-M, M))$, for some~$\alpha > 0$. In particular,~$v_{[-M, M]}$ is a pointwise solution of~\eqref{v-MMDir}. Now, we claim that
\begin{equation} \label{vI<1}
|v_{[-M, M]}(x)| < 1 \quad \mbox{for any } x \in (-M, M).
\end{equation}
Indeed, if~\eqref{vI<1} does not hold, then there exists a point~$x_0 \in (-M, M)$ at which, say,~$v_{[-M, M]}(x_0) = -1$. But then, in view of Proposition~\ref{SCPprop}, we deduce that~$v_{[-M, M]} = -1$ in the whole of~$(-M, M)$, which clearly contradicts the continuity of~$v_{[-M, M]}$ at~$M$. Thus,~\eqref{vI<1} is valid.

Given~$\tau \ge 0$, set~$u(x) := v_{[-M, M]}(x)$ and~$u_\tau(x) := v_{[-M, M]}(x - \tau)$, for any~$x \in \R$. Note that
\begin{equation} \label{utau=u}
u_\tau(x) = u(x) \quad \mbox{for any } x \in (-\infty, -M] \cup [M + \tau, +\infty).
\end{equation}
We define
$$
\hat{\tau}_0 := \inf \Big\{ \tau_0 > 0 : u_\tau \le u \mbox{ in } \R, \mbox{ for any } \tau \ge \tau_0 \Big\}.
$$
By construction, it holds~$\hat{\tau}_0 \in [0, 2M]$. Observe that if we show that
\begin{equation} \label{hattau0=0}
\hat{\tau}_0 = 0,
\end{equation}
then the monotonicity of~$v_{[-M, M]}$ would follow. To prove~\eqref{hattau0=0}, we argue by contradiction and in fact suppose that~$\hat{\tau}_0 \in (0, 2M]$. As a result,
\begin{equation} \label{uhattau0leu}
u_{\hat{\tau}_0} \le u \quad \mbox{in } \R,
\end{equation}
and there exist two sequences~$\varepsilon_j > 0$ and~$x_j \in \R$ such that~$\varepsilon_j \rightarrow 0$ as~$j \rightarrow +\infty$ and
\begin{equation} \label{uhattau0>uxj}
u_{\hat{\tau}_0 - \varepsilon_j}(x_j) > u(x_j),
\end{equation}
for any~$j \in \N$. Moreover, by~\eqref{utau=u}, we have that~$x_j \in (-M, M + \hat{\tau}_0 - \varepsilon_j)$, so that~$x_j$ converges to some~$x_0 \in [-M, M + \hat{\tau}_0]$, up to subsequences. Using~\eqref{uhattau0leu} and~\eqref{uhattau0>uxj}, it then follows that
\begin{equation} \label{uhattau0=ux0}
u_{\hat{\tau}_0}(x_0) = u(x_0),
\end{equation}
while by~\eqref{vI<1} we further deduce that~$\hat{\tau}_0 < 2M$ and~$x_0 \in (-M + \hat{\tau}_0, M)$. By virtue of~\eqref{v-MMDir},~\eqref{uhattau0leu} and~\eqref{uhattau0=ux0}, we may now apply Proposition~\ref{SCPprop} and obtain that~$u_{\hat{\tau}_0}(x) = u(x)$, for any~$x \in (-M + \hat{\tau}_0, M)$. By~\eqref{vI<1} and the continuity of~$v_{[-M, M]}$, we are then led to
$$
1 > v_{[-M, M]}(M - \hat{\tau}_0) = u_{\hat{\tau}_0}(M) = u(M) = v_{[-M, M]}(M) = 1,
$$
which is a contradiction. Accordingly,~\eqref{hattau0=0} is true and therefore~$v_{[-M, M]}$ is non-decreasing.

Now, we show that~$v_{[-M, M]}$ is the unique solution of the Dirichlet problem~\eqref{v-MMDir}. Let~$w$ be a solution of~\eqref{v-MMDir}. By Proposition~\ref{C2salphaintsemiprop}, we know that~$w \in C^\alpha(\R) \cap C^{2 s + \alpha}_\loc((-M, M))$, for some~$\alpha > 0$. Furthermore, by arguing as in the proof of~\eqref{vI<1}, we get that~$|w(x)| < 1$, for any~$x \in (-M, M)$. We claim that
\begin{equation} \label{wlev-MM}
w \le v_{[-M, M]} \quad \mbox{in } \R.
\end{equation}
To prove it, we take any~$\tau \ge 0$ and set~$w_\tau(x) := w(x - \tau)$, for any~$x \in \R$. Note that~$w_\tau(x) = v_{[-M, M]}(x)$, for any~$x \in (-\infty, -M] \cup [M + \tau, +\infty)$. Set then
$$
\bar{\tau}_0 := \inf \Big\{ \tau_0 > 0 : w_\tau \le v_{[-M, M]} \mbox{ in } \R, \mbox{ for any } \tau \ge \tau_0 \Big\} \in [0, 2 M).
$$
Clearly,~\eqref{wlev-MM} would follow if we prove that~$\bar{\tau}_0 = 0$. We thus argue by contradiction and suppose that~$\bar{\tau}_0 > 0$. Then, it is not hard to show that~$w_{\bar{\tau}} \le v_{[-M, M]}$ in~$\R$ and that there exists a point~$x_0 \in (-M + \hat{\tau}_0, M)$ at which~$w_{\bar{\tau}_0}(x_0) = v_{[-M, M]}(x_0)$. But then, by Proposition~\eqref{SCPprop} we deduce that~$w_{\bar{\tau}} = v_{[-M, M]}$ in the whole interval~$[-M + \hat{\tau}_0, M]$, which is a contradiction, since~$\hat{\tau}_0 > 0$. Accordingly,~\eqref{wlev-MM} is valid. With a completely analogous argument we obtain that the converse inequality is also true and, therefore, that~$w = v_{[-M, M]}$.

Finally, we are left to prove that~$v_{[-M, M]}$ is an odd function. To do this, we define
$$
z(x) := - v_{[-M, M]}(-x) \quad \mbox{for any } x \in \R.
$$
Clearly, we have that~$z(x) = -1$ for any~$x \le -M$ and~$z(x) = 1$ for any~$x \ge M$. Moreover,
$$
L_K z(x) = - L_K v_{[-M, M]}(-x) = - W'(v_{[-M, M]}(-x)) = - W'(-z(x)),
$$
for any~$x \in (-M, M)$. By taking advantage of~\eqref{Weven}, we have that~$W'$ is odd in~$[-1, 1]$ and we conclude that~$z$ is a solution of~\eqref{v-MMDir}. Hence,~$z = v_{[-M, M]}$, by uniqueness, and~$v_{[-M, M]}$ is odd.
\ndim

\subsection{Minimizers on the real line} \label{lineminsubsec}

We now use the results obtained in the previous subsection to deduce the existence of a class~A minimizer for~$\E_K$ in~$\R$.

Recalling definitions~\eqref{ics} and~\eqref{G*def}, we introduce the set of monotone minimizers
$$
	\mathcal{M}: = \Big\{ u \in \mathcal{X} : u \mbox{ is a non-decreasing class~A minimizer for } \E_K \Big\}.
$$

In the next proposition we show that the class~$\mathcal{M}$ defined above contains at least one element.

\bpro \label{Mnonemptyprop}
The set~$\mathcal{M}$ is not empty. In particular, there exists an odd class~A minimizer~$u_0: \R \to (-1, 1)$ for~$\E_K$, which is~$C^{1 + 2 s + \alpha}(\R^N)$ regular, for some~$\alpha > 0$, and satisfies
\begin{align}
\label{u*(0)=0} u_0(0) & = 0, \\
\label{u*'>0} u_0'(x) & > 0 \quad \mbox{for any } x \in \R, \\
\label{u*lim} \lim_{x \rightarrow \pm \infty} u_0(x) & = \pm 1,
\end{align}
and~\eqref{N=1enest}.
\npro
\bdim
Let~$M > 5$ and consider the local minimizer~$v_{[-M, M]}: \R \to [-1, 1]$ given by Lemma~\ref{intminlem}. Recall that~$v_{[-M, M]}$ is an odd, non-decreasing function such that $v_{[-M, M]}(x) = -1$ if $x \leq -M$ and $v_{[-M, M]}(x) = 1$ if $x \ge M$. Moreover,
\begin{equation} \label{vconh1}
	\frac{1}{2} \left[ v_{[-M, M]} \right]_{\mathbb{H}^K(J)}^2 \le \E_K(v_{[-M, M]}, J) \le C_1 \Psi_s(|J|),
\end{equation}
where either~$J = [-M, M]$ or~$J$ is any subinterval of~$[- M, M]$, with~$|J| > 6$ and~$\dist(J, \R \setminus [-M, M]) > 2$. Note that~$C_1 \ge 1$ is a constant depending only on~$s$,~$\Lambda$ and~$W$. Also,~$v_{[-M, M]}$ is a solution of
\begin{equation} \label{vconh2}
	L_K v_{[-M, M]} = W'(v_{[-M, M]}) \quad \mbox{in } (-M, M),
\end{equation}
and thus by Proposition~\ref{C2salphaintsemiprop} we deduce that~$v_{[-M, M]} \in C^\alpha(\R)$, for some~$\alpha \in (0, 1)$, with H\"{o}lder norm bounded independently\footnote{A careful inspection of the proof of~\cite[Proposition~1.1]{R-OS14} - on which Propositions~\ref{C2salphaintsemiprop} is based - shows that the H\"{o}lder norm of the solution of the Dirichlet problem~\eqref{Dirichlet problem} is bounded by a constant that does not depend on~$\Omega$ as a whole, but only on the~$C^{1, 1}$ norm of its boundary (see also~\cite[Subsection~3.4.1]{C16}). In particular, when~$N = 1$ the constant is independent on the reference interval. As a result, we can conclude that the~$C^\alpha(\R)$ norm of~$v_{[-M, M]}$ is independent of~$M$.} of~$M$.

In view of this and Ascoli-Arzel\`a's theorem, we may assume that~$v_{[-M, M]}$ converges to a continuous function~$u_0$, uniformly on compacts subsets of~$\R$, as~$M \rightarrow +\infty$. By the oddness~$v_{[-M, M]}$, we have that~$v_{[-M, M]}(0) = 0$, for any~$M$. Accordingly,~$u_0$ satisfies~\eqref{u*(0)=0}. Also,~$u_0$ is odd, non-decreasing and weakly satisfies
\begin{equation} \label{u*eq}
L_K u_0 = W'(u_0) \quad \mbox{in } \R, 
\end{equation}
in view of~\eqref{vconh1},~\eqref{vconh2} and Lemma~\ref{stablem}.  By Proposition~\ref{entiresemiregprop}, it then follows that~$u_0 \in C^{1 + 2 s + \alpha}(\R)$, for some~$\alpha > 0$.

Now we prove that~$u_0 \in \mathcal{M}$, thus concluding the proof of the proposition. In order to do this, we first show that~\eqref{N=1enest} holds true. To check it, we fix~$R > 4$ and address the energy of~$v_{[-M, M]}$ inside the interval~$[-R, R]$. By taking~$M$ suitably large in dependence of~$R$ if necessary, by~\eqref{vconh1} we have that
$$
\E_K(v_{[-M, M]}, [-R, R]) \le C \Psi_s(R),
$$
for some constant~$C > 0$ independent of~$M$ and~$R$. The finiteness condition~\eqref{N=1enest} then follows by letting~$R$ go to~$+\infty$ in the above inequality, thanks to Fatou's lemma.

Next, we check that~\eqref{u*lim} holds true. In view of the monotonicity of~$u_0$ and~\eqref{u*(0)=0}, we know that there exist two numbers~$-1 \le a_- \le 0 \le a_+ \le 1$ such that
$$
\lim_{x \rightarrow \pm \infty} u_0(x) = a_\pm.
$$
We prove here that~$a_+ = 1$, while a completely analogous argument shows that~$a_- = -1$ holds too. Suppose by contradiction that~$a_+ < 1$ and notice that~$u_0(x) \in [0, a_+]$ for any~$x \ge 0$. Set
$$
\kappa := \inf_{x \ge 0} W(u_0(x)) = \inf_{r \in [0, a_+]} W(r).
$$
By taking advantage of~\eqref{Wpos} in combination with the fact that~$a_+ < 1$, we deduce that~$\kappa > 0$. Consequently,
$$
\G^*(u_0) \ge \limsup_{R \rightarrow +\infty} \frac{1}{\Psi_s(R)} \int_{0}^{R} W(u_0(x)) \, dx \ge \kappa \lim_{R \rightarrow +\infty} \frac{R}{\Psi_s(R)} = +\infty,
$$
in contradiction with~\eqref{N=1enest}. Thence,~\eqref{u*lim} is valid. In particular,~$u_0 \in \mathcal{X}$.

Finally, the monotonicity of~$u_0$,~\eqref{u*eq},~\eqref{u*lim} and Lemma~\ref{rk1pflem} imply that~$u_0$ satisfies~\eqref{u*'>0}. By virtue of this,~\eqref{u*eq} and~\eqref{u*lim}, the function~$u_0$ fulfills the hypotheses of Theorem~\ref{theorem3}. Therefore, it follows that~$u_0$ is a class~A minimizer. By this and again~\eqref{u*'>0}, we conclude that~$u_0 \in \mathcal{M}$. The proof of the proposition is thus complete.
\ndim

Next, we address the problem of assessing how~\emph{big} the set~$\mathcal{M}$ is. In Proposition~\ref{Mnonemptyprop}, we have established that~$\mathcal{M}$ contains at least one element~$u_0$. Clearly, it also contains the translations~$u_0(\cdot - k)$, for any~$k \in \R$. We are thence led to study the subclasses
$$
\mathcal{M}_{x_0} := \Big\{ u \in \mathcal{M} \, : \, x_0 = \sup \, \{ x \in \R : u(x) < 0 \} \Big\},
$$
for any fixed~$x_0 \in \R$. Of course, we have that
$$
\mathcal{M} = \bigcup_{x_0 \in \R} \mathcal{M}_{x_0} \quad \mbox{and} \quad \mathcal{M}_{x_0} \cap \mathcal{M}_{x_1} = \varnothing, \mbox{ if } x_0 \ne x_1.
$$
Also,
\begin{equation} \label{Mx0rel}
u \in \mathcal{M}_{x_0} \quad \mbox{if and only if} \quad u(\cdot + x_0) \in \mathcal{M}_{0},
\end{equation}
for any~$x_0 \in \R$. It turns out that each of these subclasses is a singleton, as shown by the following

\bpro
For any fixed~$x_0 \in \R$, the class~$\mathcal{M}_{x_0}$ consists of one single element~$u_{x_0}$. More specifically,~$u_{x_0}: \R \to (-1, 1)$ is a class~A minimizer for~$\E_K$, which is~$C^{1 + 2 s + \alpha}$ regular, for some~$\alpha > 0$, and satisfies~\eqref{u*'>0},~\eqref{u*lim},~\eqref{N=1enest} and~$u_{x_0}(x_0) = 0$.
\npro
\bdim
In light of~\eqref{Mx0rel}, it is enough to prove the statement for the point~$x_0 = 0$. Note that the function~$u_0$ constructed in Proposition~\ref{Mnonemptyprop} belongs to~$\mathcal{M}_{0}$. Let~$u \in \mathcal{M}_{0}$. Observe that~$u$ is a weak solution of~\eqref{u*eq} and, hence, by Proposition~\ref{entiresemiregprop}, that~$u \in C^{1 + 2 s + \alpha}(\R)$, for some~$\alpha > 0$. Also,~$|u| \le 1$ in~$\R$, since~$u$ satisfies~\eqref{u*lim} and it is non-decreasing. If we show that~$u = u_0$, then the proof would be over.

First, we notice that there exists a small value~$\varepsilon_0 > 0$ such that, for any~$\varepsilon \in (0, \varepsilon_0)$, we can find a~$\bar{k}_\varepsilon \in \R$ for which
\begin{equation} \label{startslideuu*}
\mbox{if } k \le \bar{k}_\varepsilon, \mbox{ then } u(x - k) + \varepsilon > u_0(x) \mbox{ for any } x \in \R.
\end{equation}
This is true as a consequence of both~$u_0$ and~$u$ having values in~$[-1, 1]$ and satisfying~\eqref{u*lim}. Then, we start sliding the graph of~$u + \varepsilon$ to the right until it first touches that of~$u_0$. That is, we take~$\hat{k}_\varepsilon$ as the largest possible value of~$\bar{k}_\varepsilon$ for which~\eqref{startslideuu*} holds true, and find a point~$x_\varepsilon \in \R$ at which
$$
u(x_\varepsilon - \hat{k}_\varepsilon) + \varepsilon = u_0(x_\varepsilon).
$$
Again, this is possible in view of the continuity and the behaviour at~$\pm \infty$ of~$u$ and~$u_0$. Set now~$u_\varepsilon(x) := u(x - \hat{k}_\varepsilon) + \varepsilon$ and observe that, by definition of~$\hat{k}_\varepsilon$, it holds
\begin{equation} \label{uepsproperties}
\begin{cases}
u_\varepsilon(x) \ge u_0(x) \quad \mbox{for any } x \in \R \\
u_\varepsilon(x_\varepsilon) = u_0(x_\varepsilon).
\end{cases}
\end{equation}

Now we claim that
\begin{equation} \label{xepsbounded}
x_\varepsilon \mbox{ is bounded as } \varepsilon \rightarrow 0^+.
\end{equation}
By contradiction, suppose that there is a sequence of values~$\varepsilon_j > 0$ for which~$\varepsilon_j \rightarrow 0^+$ and, say,~$x_{\varepsilon_j} \rightarrow +\infty$, as~$j \rightarrow +\infty$. By~\eqref{W''>0}, we can pick a small value~$c > 0$ such that~$W'$ is monotone non-decreasing in~$[1 - c, 1]$. Fix a real number~$M > 0$ large enough to have~$u_0(M) > 1 - c/2$. Notice that~$x_{\varepsilon_j} > M$ and~$\varepsilon_j < c/2$, provided~$j$ is sufficiently large. By this and the monotonicity of~$u_0$, we have
$$
u_0(x) \ge u_0(x) - \varepsilon_j \ge u_0(M) - \frac{c}{2} > 1 - c,
$$
for any~$x \in (M, +\infty)$. Hence, recalling the monotonicity of~$W'$ in~$[1 - c ,1]$, we obtain
\begin{equation} \label{W'monot}
W'(u_0(x) - \varepsilon_j) \le W'(u_0(x)) \quad \mbox{for any } x \in (M, +\infty).
\end{equation}
Observe now that, since both~$u$ and~$u_0$ satisfy~\eqref{u*eq},
$$
\begin{cases}
L_K u_{\varepsilon_j} = W'(u_{\varepsilon_j} - \varepsilon_j) & \quad \mbox{in } \R \\
L_K u_0 = W'(u_0) & \quad \mbox{in } \R.
\end{cases}
$$
Consequently, by this and~\eqref{uepsproperties}, we are able to use Proposition~\ref{SCPprop} - with~$\Omega = (M, +\infty)$,~$f_1(r) = W'(r - \varepsilon_j)$ and~$f_2(r) = W'(r)$ - to deduce that
\begin{equation} \label{identcontrad}
u_{\varepsilon_j}(x) = u_0(x) \quad \mbox{for any } x > M,
\end{equation}
provided~$j$ is large enough. Notice that the validity of condition~\eqref{f1gef2} there is ensured by~\eqref{W'monot}. But~\eqref{identcontrad} is contradictory, as can be seen for instance by letting~$x \rightarrow +\infty$. A symmetrical argument shows that we reach a contradiction also if~$x_{\varepsilon_j} \rightarrow -\infty$. Thus,~\eqref{xepsbounded} follows. As a result of~\eqref{xepsbounded}, we have that, up to a subsequence,
\begin{equation} \label{xepslimit}
\lim_{\varepsilon \rightarrow 0^+} x_\varepsilon = x_0,
\end{equation}
for some~$x_0 \in \R$.

Then, we claim that
\begin{equation} \label{kepsbounded}
\hat{k}_\varepsilon \mbox{ is bounded as } \varepsilon \rightarrow 0^+.
\end{equation}
Again, we argue by contradiction and suppose that~$\hat{k}_{\varepsilon_j} \rightarrow \pm \infty$ on an infinitesimal sequence~$\varepsilon_j > 0$. Applying the identity on the second line of~\eqref{uepsproperties} and~\eqref{xepslimit}, we obtain
$$
\mp 1 = \lim_{j \rightarrow +\infty} \left[ u(x_{\varepsilon_j} - \hat{k}_{\varepsilon_j}) + \varepsilon_j \right] = \lim_{j \rightarrow +\infty} u_{\varepsilon_j}(x_{\varepsilon_j}) = \lim_{j \rightarrow +\infty} u_0(x_{\varepsilon_j}) = u_0(x_0),
$$
which is not the case, since~$u_0$ has values in~$(-1, 1)$. Thence,~\eqref{kepsbounded} holds and, up to a subsequence,
\begin{equation} \label{kepslimit}
\lim_{\varepsilon \rightarrow 0^+} \hat{k}_\varepsilon = \hat{k}_0,
\end{equation}
for some~$\hat{k}_0 \in \R$.

By virtue of~\eqref{xepslimit} and~\eqref{kepslimit}, we may finally let~$\varepsilon \rightarrow 0^+$ in~\eqref{uepsproperties}, to find that
$$
\begin{cases}
L_K u(x - \hat{k}_0) = W'(u(x - \hat{k}_0)) & \quad \mbox{for any } x \in \R \\
L_K u_0(x) = W'(u_0(x)) & \quad \mbox{for any } x \in \R \\
u(x - \hat{k}_0) \ge u_0(x) & \quad \mbox{for any } x \in \R \\
u(x_0 - \hat{k}_0) = u_0(x_0). &
\end{cases}
$$
By applying once again Proposition~\ref{SCPprop}, we infer that~$u(x - \hat{k}_0) = u_0(x)$ for any~$x \in \R$. Then, as~$u, u_0 \in \mathcal{M}_{0}$, we conclude that~$u(0) = 0 = u_0(0)$. Recalling~\eqref{u*'>0}, it follows that~$\hat{k}_0 = 0$ and hence~$u = u_0$. The proposition is thus proved.
\ndim

\subsection{Further estimates: general kernels} \label{est1subsec}

Up to now, we have established the existence - and essential uniqueness - of the minimizer~$u_0$ in the class~$\mathcal{X}$. Moreover, we already know by construction that~$u_0$ is strictly increasing and that~\eqref{N=1enest} holds true.

In this subsection we show that estimates~\eqref{stima1} and~\eqref{N=1tailest} are also valid. These results are the content of the following two propositions.

\bpro
The function~$u_0$ constructed in Proposition~\ref{Mnonemptyprop} satisfies the decay estimates~\eqref{stima1}.
\npro
\bdim
We begin by addressing the validity of the first estimate in~\eqref{stima1}. Obviously, we may restrict ourselves to prove only that there exists~$R_1, C_1 > 0$ such that
\begin{equation} \label{udecay-infty}
u_0(x) \le - 1 + \frac{C_1}{|x|^{2 s}},
\end{equation}
if~$x \le - R_1$.

To do this, first observe that, by~\eqref{W''>0},
\begin{equation} \label{W'geW'}
W'(t) \ge W'(r) + c (t - r) \quad \mbox{for any } r \le t \mbox{ such that } r, t \in [-1, -1 + c],
\end{equation}
for some~$c \in (0, 1/2)$. Take now~$\tau = c$ in Lemma~\ref{barrierlem} and for any~$R \ge C$ consider the barrier~$w$ constructed there. By~\eqref{wbarrange},~\eqref{wbar=1} and~\eqref{u*lim}, there exists~$k_0 \in \R$ such that
\begin{equation} \label{u0<wbar}
\mbox{for any } k \in (-\infty, k_0), \mbox{ it holds } u_0(x) < w(x - k) \mbox{ for any } x \in \R.
\end{equation}
Now, let~$\bar{k}_0$ be the largest~$k_0$ for which~\eqref{u0<wbar} is true. Clearly,
\begin{equation} \label{u0lewbar}
u_0(x) \le w(x - \bar{k}_0) \quad \mbox{for any } x \in \R.
\end{equation}
Also, it is not hard to check that there exists
\begin{equation} \label{barxdomain}
\bar{x} \in (\bar{k}_0 - R, \bar{k}_0 + R),
\end{equation}
at which
\begin{equation} \label{u0=wbar}
u_0(\bar{x}) = w(\bar{x} - \bar{k}_0).
\end{equation}

We claim that
\begin{equation} \label{u0barxge-1+c}
u_0(\bar{x}) \ge - 1 + c.
\end{equation}
To prove it, we argue by contradiction and suppose indeed that
\begin{equation} \label{u0barx<-1+x}
u_0(\bar{x}) \in (-1, -1 + c).
\end{equation}
Define
$$
\Omega := \left\{ x \in (\bar{k}_0 - R, \bar{k}_0 + R) : u_0(x) < -1 + c \right\},
$$
and note that, by~\eqref{u0barx<-1+x} and the continuity and monotonicity of~$u_0$, we have that~$\Omega$ is an open domain with
\begin{equation} \label{Omegainclusion}
(\bar{k}_0 - R, \bar{x}] \subset \Omega.
\end{equation}
Setting now~$\bar{w}(x) := w(x - \bar{k}_0)$, by~\eqref{LKwbar},~\eqref{u*eq},~\eqref{u0lewbar} and~\eqref{u0=wbar}, we know that
$$
\begin{cases}
L_K \bar{w} \le c(1 + \bar{w}) & \quad \mbox{in } (\bar{k}_0 - R, \bar{k}_0 + R) \\
L_K u_0 = W'(u_0) & \quad \mbox{in } \R \\
\bar{w} \ge u_0 & \quad \mbox{in } \R \\
\bar{w}(\bar{x}) = u_0(\bar{x}).
\end{cases}
$$
Furthermore, notice that, by taking~$t = u_0(x)$ and~$r = -1$ in~\eqref{W'geW'} and recalling~\eqref{W'=0},
$$
W'(u_0(x)) \ge c (1 + u_0(x)) \quad \mbox{for any } x \in \Omega.
$$
In view of this last consideration, we are then in position to apply Proposition~\ref{SCPprop} and obtain that~$u_0(x) = \bar{w}(x)$, for any~$x \in \Omega$. But then, by~\eqref{Omegainclusion}, the continuity of~$u_0, \bar{w}$ and~\eqref{wbar=1},
$$
1 > u_0(\bar{k}_0 - R) = \bar{w}(\bar{k}_0 - R) = w(-R) = 1,
$$
which is a contradiction. Consequently,~\eqref{u0barxge-1+c} holds true.

In view of~\eqref{1+wbarest},~\eqref{u0=wbar},~\eqref{barxdomain} and~\eqref{u0barxge-1+c} we now get
$$
C(R + 1 - |\bar{x} - \bar{k}_0|)^{-2 s} \ge 1 + w(\bar{x} - \bar{k}_0) = 1 + u_0(\bar{x}) \ge c,
$$
so that
\begin{equation} \label{barx-bark0ge}
|\bar{x} - \bar{k}_0| \ge R - c',
\end{equation}
for some~$c' > 0$. Moreover,
\begin{equation} \label{barxgebark0}
\bar{x} \ge \bar{k}_0.
\end{equation}

To check~\eqref{barxgebark0}, we argue by contradiction and suppose that~$\bar{x} < \bar{k}_0$. Set~$\hat{k} := 2 \bar{x} - \bar{k}_0$ and notice then that~$\hat{k} < \bar{k}_0$. Accordingly, by~\eqref{u0<wbar} and~\eqref{u0=wbar} we deduce that
$$
w(\hat{k}_0 - \bar{x}) = w(\bar{x} - \hat{k}) > u_0(\bar{x}) = w(\bar{x} - \hat{k}_0),
$$
in contradiction with the parity of~$w$. Thus,~\eqref{barxgebark0} is true.

In consequence of~\eqref{barxdomain},~\eqref{barx-bark0ge} and~\eqref{barxgebark0}, we see that
\begin{equation} \label{barx-bark0int}
\bar{x} - \bar{k}_0 \in [R - c', R].
\end{equation}
Let~$\kappa > 0$ be chosen in such a way that~$u_0(-\kappa) = -1 + c$. By the monotonicity of~$u_0$, we clearly have~$- \kappa \le \bar{x}$ and
\begin{equation} \label{u0kappaleu0barx}
u_0(x - \kappa) \le u_0(x + \bar{x}) \quad \mbox{for any } x \in \R.
\end{equation}
Take now any~$y \in [R/2, R]$. By~\eqref{barx-bark0int} and taking a larger~$R$ if necessary, we have that~$\bar{x} - y - \bar{k}_0 \in [-R/2, R/2]$. Consequently, by~\eqref{1+wbarest},
$$
1 + w(\bar{x} - y - \bar{k}_0) \le C (R + 1 - |\bar{x} - y - \bar{k}_0|)^{-2 s} \le C \left( \frac{R}{2} \right)^{-2 s} \le 4 C y^{-2 s}.
$$
By combining this with~\eqref{u0lewbar} and~\eqref{u0kappaleu0barx}, we then get
$$
u_0(- \kappa - y) \le u_0(\bar{x} - y) \le w(\bar{x} - y - \bar{k}_0) \le - 1 + 4 C y^{- 2 s} \quad \mbox{for any } y \in \left[ \frac{R}{2}, R \right].
$$
Since~$\kappa$ is a positive constant and~$R$ may be chosen arbitrarily large, it is almost immediate to check that this implies~\eqref{udecay-infty}. Accordingly, the first estimate in~\eqref{stima1} is established.

Now, we head to the proof of the second estimate of~\eqref{stima1}. We first remark that, since~$u_0 \in C^{1 + 2 s + \alpha}(\R)$, for some~$\alpha > 0$, we may differentiate equation~\eqref{u*eq} and deduce that~$u_0'$ solves
\begin{equation} \label{Lu'}
L_K u_0' = W''(u_0) u_0' \quad \mbox{in } \R.
\end{equation}
Observe now that, in view of~\eqref{W''>0} and the fact that~$u_0$ satisfies~\eqref{u*lim}, we can take~$R_0 > 0$ big enough to have~$W''(u_0(x)) \ge \delta$ for any~$x \in \R$ such that~$|x| > R_0$ and for some constant~$\delta > 0$. By virtue of this,~\eqref{Lu'} and~\eqref{u*'>0}, we then obtain that
$$
L_K u_0' \ge \delta u_0' \quad \mbox{in } \R \setminus [-R_0, R_0].
$$
The thesis now follows from Lemma~\ref{vsubdecay}.
\ndim

\bpro \label{UTEprop}
The upper tail energy estimate~\eqref{N=1tailest} holds true.
\npro
\bdim
All along the proof, we denote with~$c$ any positive constant, whose value may change from line to line.

First we notice that, by the second estimate in~\eqref{stima1},
\begin{equation} \label{u'Linftydecay}
\| u_0' \|_{L^\infty \left( \left[ t - \frac{|t|}{2}, t + \frac{|t|}{2} \right] \right) } \le \frac{c}{|t|^{1 + 2 s}},
\end{equation}
for~$|t|$ sufficiently large. Moreover, given any~$\rho > 0$, by the fact that~$|u_0| \le 1$, we compute 
\begin{equation} \label{rhoenest}
\begin{aligned}
\int_{\R} \frac{\left| u_0(x) - u_0(t) \right|^2}{|x - t|^{1 + 2 s}} \, dx & \le 2 \left( \| u_0' \|_{L^\infty([t - \rho, t + \rho])}^2 \int_0^{t + \rho} \frac{dx}{|x - t|^{- 1 + 2 s}} + \int_{t + \rho}^{+\infty} \frac{4 \, dx}{|x - t|^{1 + 2 s}} \right) \\
& \le c \left( \| u_0' \|_{L^\infty([t - \rho, t + \rho])}^2 \rho^2 + 1 \right) \rho^{- 2 s}.
\end{aligned}
\end{equation}

We claim that
\begin{equation} \label{betaest}
\beta(t) := \frac{1}{4} \int_{\R} \left| u_0(x) - u_0(t) \right|^2 K(x - t) \, dx + W(u_0(t)) \le \frac{c}{1 + |t|^{2 s}},
\end{equation}
for any~$t \in \R$. We actually prove the stronger
\begin{equation} \label{beta2est}
\frac{1}{4} \int_{\R} \frac{\left| u_0(x) - u_0(t) \right|^2}{|x - t|^{1 + 2 s}} \, dx + W(u_0(t)) \le \frac{c}{1 + |t|^{2 s}},
\end{equation}
for any~$t \in \R$. Observe that~\eqref{beta2est} implies~\eqref{betaest}, thanks to the right-hand inequality of~\eqref{weakellipticity}.

To prove~\eqref{beta2est}, we first plug~$\rho = |t| / 2$ into~\eqref{rhoenest}. In view of~\eqref{u'Linftydecay} we get
\begin{equation} \label{kinendecay}
\int_{\R} \frac{\left| u_0(x) - u_0(t) \right|^2}{|x - t|^{1 + 2 s}} \, dx \le \frac{c}{|t|^{2 s}},
\end{equation}
provided~$|t|$ is large enough. Also,~$u_0' \in L^\infty(\R)$ and thus, by choosing e.g.~$\rho = 1$ in~\eqref{rhoenest},
\begin{equation} \label{kinenbound}
\int_{\R} \frac{\left| u_0(x) - u_0(t) \right|^2}{|x - t|^{1 + 2 s}} \, dx \le c,
\end{equation} 
for any~$t \in \R$. On the other hand,~$W$ is of class~$C^2$ and satisfies~\eqref{W'=0}. Hence, recalling the first estimate of~\eqref{stima1} we obtain
\begin{align*}
W(u_0(t)) & = W(u_0(t)) - W(1) = \int_1^{u_0(t)} W'(\tau) \, d\tau = \int_{u_0(t)}^1 [W'(1) - W'(\tau)] \, d\tau \\
& \le \| W'' \|_{L^\infty([-1, 1])} \int_{u_0(t)}^1 (1 - \tau) \, d\tau = \frac{\| W'' \|_{L^\infty([-1, 1])}}{2} (1 - u_0(t))^2 \\
& \le \frac{c}{|t|^{4 s}},
\end{align*}
if~$t$ is close enough to~$1$. Similarly, one prove that the same is true when~$t$ approaches~$-1$. By this and the boundedness of~$W$ we get that
\begin{equation} \label{W(u0)est}
W(u_0(t)) \le \frac{c}{1 + |t|^{4 s}},
\end{equation}
for any~$t \in \R$. The combination of~\eqref{kinendecay},~\eqref{kinenbound} and~\eqref{W(u0)est}  leads to~\eqref{beta2est}.

With the aid of the previous computations, we may now head to the actual proof of~\eqref{N=1tailest}. We have
\begin{align*}
\int_{-R}^{R} \int_{\R \setminus [-R, R]} \left| u_0(x) - u_0(y) \right|^2 K(x - y) \, dx dy & \le \int_{-\frac{R}{2}}^{\frac{R}{2}} \left( \int_{R}^{+\infty} \frac{8 \Lambda \, dx}{|x - y|^{1 + 2 s}} \right) dy + 4 \int_{ \left\{ \frac{R}{2} < |y| \le R \right\} } \beta(y) \, dy \\
& \le c \left[ \int_{-\frac{R}{2}}^{\frac{R}{2}} (R - y)^{-2 s} \, dy + \int_{\frac{R}{2}}^{R} \frac{dy}{1 + y^{2 s}} \right] \\
& \le c R^{1 - 2 s}.
\end{align*}
This finishes the proof of the proposition.
\ndim

Notice that we did not really need inequality~\eqref{W(u0)est} to prove Proposition~\ref{UTEprop}. However, we included such estimate for the potential term, as it will turn out to be helpful later in Section~\ref{thm2sec}.

\subsection{Further estimates: positive kernels} \label{est2subsec}
Here we tackle~\eqref{N=1tailestbelow} and~\eqref{N=1enestbelow}. Since both of them are estimates from below, to prove them we assume the more restrictive condition~\eqref{ellipticity} on~$K$. Thus,~\eqref{ellipticity} will be implicitly required throughout the subsection.

\bpro
The lower tail energy estimate~\eqref{N=1tailestbelow} holds true.
\npro
\bdim
Let~$R > 0$ be large enough to have
$$
u_0(x) \ge \frac{1}{2} \mbox{ for any } x \ge R \qquad \mbox{and} \qquad  u_0(y) \le - \frac{1}{2} \mbox{ for any } y \le - \frac{R}{4}.
$$
For such values of~$R$, using~\eqref{ellipticity} we compute
\begin{equation} \label{N=1tailtech}
\begin{aligned}
\int_{-\frac{R}{2}}^{-\frac{R}{4}} \int_{R}^{+\infty} \left| u_0(x) - u_0(y) \right|^2 K(x - y) \, dx dy & \ge \lambda \int_{-\frac{R}{2}}^{-\frac{R}{4}} \left( \int_{R}^{+\infty} \frac{dx}{|x - y|^{1 + 2 s}} \right) dy \\
& = \frac{\lambda}{2 s} \int_{-\frac{R}{2}}^{-\frac{R}{4}} (R - y)^{- 2 s} dy \\
& \ge \frac{\lambda}{2^{3 + 2 s} s} R^{1 - 2 s}.
\end{aligned}
\end{equation}
Formula~\eqref{N=1tailest} then immediately follows.
\ndim

We conclude this subsection with a lemma that gives a sharp lower bound for the total energy~$\E_K(u_0, [-R, R])$, when~$s = 1/2$.

\blem \label{s=1/2belowlem}
Let~$s = 1/2$. There exists a constant~$c > 0$ such that
\begin{equation} \label{s=1/2belowest}
\E_K(u_0, [-R, R]) \ge c \log R,
\end{equation}
for any~$R$ large enough.
\nlem
\bdim
Choose~$k_0 > 1$ in a way that
\begin{equation} \label{k0def}
u_0(x) \ge \frac{1}{2} \mbox{ for any } x \ge k_0 \qquad \mbox{and} \qquad u_0(y) \le - \frac{1}{2} \mbox{ for any } y \le - k_0.
\end{equation}
Let~$\ell > k \ge k_0$ and define
$$
I_{k, \ell} := \int_{-\ell}^{-k} \int_k^\ell \left| u_0(x) - u_0(y) \right|^2 K(x - y) \, dx dy.
$$
By~\eqref{k0def} and~\eqref{ellipticity} we compute
$$
I_{k, \ell} \ge \lambda \int_{-\ell}^{-k} \left( \int_k^{\ell} \frac{dx}{(x - y)^{2}} \right) dy = \lambda \int_{-\ell}^{-k} \left( \frac{1}{k - y} - \frac{1}{\ell - y} \right) dy = \lambda \log \frac{(k + \ell)^2}{4 k \ell}.
$$
If we set~$\ell = 10 k$, the above inequality becomes
\begin{equation} \label{k10kest}
I_{k, 10 k} \ge \lambda \log \frac{121 k^2}{40 k^2} > \lambda.
\end{equation}
Take now any~$R$ satisfying
\begin{equation} \label{Rboundk0}
R > 100 k_0^2,
\end{equation}
and let~$M > 0$ be the largest integer for which~$10^M k_0 \le R$. Notice that then
$$
10^{M + 1} k_0 > R,
$$
which, along with~\eqref{Rboundk0}, implies
$$
M > \log_{10} \frac{R}{k_0} - 1 = \frac{\log \frac{R}{k_0} - \log 10}{\log 10} \ge \frac{\log R}{2 \log 10}
$$
By this and~\eqref{k10kest}, we conclude that
$$
\int_{-R}^R \int_{-R}^{R} \left| u_0(x) - u_0(y) \right|^2 K(x - y) \, dx dy \ge I_{k_0, 10^M k_0} \ge \sum_{j = 1}^M I_{10^{j - 1} k_0, 10^j k_0} \ge \lambda M \ge \frac{\lambda}{2 \log 10} \log R,
$$
which gives~\eqref{s=1/2belowest}.
\ndim

Notice that we can now conclude that~\eqref{N=1enestbelow} is true. Indeed, when~$s > 1/2$ this is obvious (see Remark~\ref{s=1/2limexrk}). On the other hand, if~$s < 1/2$ this fact immediately follows from~\eqref{N=1tailestbelow}, while for~$s = 1/2$ it is a consequence of Lemma~\ref{s=1/2belowlem}.

\subsection{Further estimates: homogeneous kernels} \label{est3subsec}

Finally, we address the validity of~\eqref{barEs=1/2}. To this aim, we suppose~$s = 1/2$. Unfortunately, we are able to prove such result only for homogeneous kernels, that is - since~$N = 1$ - only for those kernels which are multiples of the kernel of the fractional Laplacian.

\bpro \label{Ebars=1/2prop}
Let~$s = 1/2$ and suppose that~$K$ is in the form~\eqref{N=1Khom}. Then,~\eqref{barEs=1/2} holds true.
\npro
\bdim
First of all, we remark that, in view of the right-hand inequality in~\eqref{N=1tailest}, we already know that
$$
\lim_{R \rightarrow +\infty} \frac{1}{\log R} \int_{-R}^R \int_{\R \setminus [-R, R]} \left| u_0(x) - u_0(y) \right|^2 K(x - y) \, dx dy = 0.
$$
Hence,
$$
\lim_{R \rightarrow +\infty} \frac{\E_K(u_0, [-R, R])}{\log R} = \lim_{R \rightarrow +\infty} \frac{\displaystyle \int_{-R}^R \beta(x) \, dx}{\log R},
$$
with~$\beta$ as in~\eqref{betaest}.

To compute this limit, we use L'H\^{o}pital's rule. Observe that we are allowed to use such method, since, by Lemma~\ref{s=1/2belowlem}, the numerator of the quotient written above diverges, as~$R \rightarrow +\infty$. By the Fundamental Theorem of Calculus, we have
\begin{equation} \label{deLH}
\lim_{R \rightarrow +\infty} \frac{\E_K(u_0, [-R, R])}{\log R} = \lim_{R \rightarrow +\infty} \frac{\displaystyle \frac{d}{dR} \int_{-R}^{R} \beta(x) \, dx}{\displaystyle \frac{d}{dR} \log R} = \lim_{R \rightarrow +\infty} R \Big( \beta(R) + \beta(-R) \Big).
\end{equation}
Now, we show that
\begin{equation} \label{betaRlim}
\lim_{R \rightarrow +\infty} R \, \beta(\pm R) = \frac{\lambda_\star}{4} \left( \lim_{x \rightarrow +\infty} u_0(x) - \lim_{x \rightarrow -\infty} u_0(x) \right)^2 = \lambda_\star.
\end{equation}
Notice that~\eqref{deLH} and~\eqref{betaRlim} immediately lead to~\eqref{barEs=1/2}.

We only deal with the limit of~$R \beta(R)$ in~\eqref{betaRlim}, the term with the minus sign being completely analogous. We claim that
\begin{equation} \label{RWvanish}
\lim_{R \rightarrow +\infty} R \, W(u_0(R)) = 0,
\end{equation}
and
\begin{equation} \label{Rkinvanish}
\lim_{R \rightarrow +\infty} R \int_{-1}^{+\infty} \left| u_0(R) - u_0(y) \right|^2 K(R - y) \, dy = 0.
\end{equation}

Observe that~\eqref{RWvanish} immediately follows from estimate~\eqref{W(u0)est}. On the other hand, to prove~\eqref{Rkinvanish}, we fix~$k_0 > 0$ large enough to have, by~\eqref{stima1},
$$
|u_0'(t)| \le \frac{c_3}{t^2} \quad \mbox{for any } t \ge k_0,
$$
for some~$c_3 > 0$. Then,
$$
\left| u_0(R) - u_0(y) \right|^2 \le \left| \int_y^R |u_0'(t)| \, dt \right|^2 \le c_3^2 \left| \int_y^R \frac{dt}{t^2} \right|^2 = c_3^2 \left| \frac{1}{y} - \frac{1}{R} \right|^2 = c_3^2 \frac{(R - y)^2}{R^2 y^2},
$$
for any~$y \ge k_0$, so that, by the right-hand inequality in~\eqref{weakellipticity},
\begin{equation} \label{kinvantech1}
\int_{k_0}^{+\infty} \left| u_0(R) - u_0(y) \right|^2 K(R - y) \, dy \le \frac{c_3^2 \, \Lambda}{R^2} \int_{k_0}^{+\infty} \frac{dy}{y^2} = \frac{c_3^2 \, \Lambda}{k_0 R^2}.
\end{equation}
Also, since~$|u_0| \le 1$, by choosing~$R > 2 k_0$ we get
\begin{equation} \label{kinvantech2}
\int_{-1}^{k_0} \left| u_0(R) - u_0(y) \right|^2 K(R - y) \, dy \le 4 \Lambda \int_{-1}^{k_0} \frac{dy}{(R - y)^2} = \frac{4 \Lambda (1 + k_0)}{(R - k_0) (R + 1)} \le \frac{8 \Lambda (1 + k_0)}{R^2}.
\end{equation}
Estimates~\eqref{kinvantech1} and~\eqref{kinvantech2} combined yield~\eqref{Rkinvanish}.

In view of~\eqref{RWvanish} and~\eqref{Rkinvanish}, we end up with
$$
\lim_{R \rightarrow +\infty} R \, \beta(R) = \frac{1}{4} \lim_{R \rightarrow +\infty} R \int_{-\infty}^{-1} \left| u_0(R) - u_0(y) \right|^2 K(R - y) \, dy.
$$
By changing variables as~$y = R (1 - z)$, this becomes
\begin{equation} \label{Rbetatech}
\lim_{R \rightarrow +\infty} R \, \beta(R) = \frac{1}{4} \lim_{R \rightarrow +\infty} R^2 \int_{1 + \frac{1}{R}}^{+\infty} \left| u_0(R) - u_0(R (1 - z)) \right|^2 K(R z) \, dz.
\end{equation}
Note that so far we never used that~$K$ is in the form~\eqref{N=1Khom}, but only the growth assumption in~\eqref{weakellipticity}. We do it now. By taking advantage of~\eqref{N=1Khom}, formula~\eqref{Rbetatech} reduces to
$$
\lim_{R \rightarrow +\infty} R \, \beta(R) = \frac{\lambda_\star}{4} \lim_{R \rightarrow +\infty} \int_{1}^{+\infty} \phi_R(z) \, dz,
$$
where
$$
\phi_R(z) := \frac{\left| u_0(R) - u_0(R (1 - z)) \right|^2}{z^2} \chi_{\left( 1 + \frac{1}{R}, +\infty \right)}(z) \quad \mbox{for a.a. } z \in (1, +\infty).
$$
Observe that
$$
\left| \phi_R(z) \right| \le \frac{4}{z^2} \in L^1((1, +\infty)),
$$
and
$$
\lim_{R \rightarrow +\infty} \phi_R(z) = \frac{\left| {\displaystyle \lim_{x \rightarrow +\infty} u_0(x) - \lim_{x \rightarrow -\infty} u_0(x)} \right|^2}{z^2} = \frac{4}{z^2},
$$
for any~$z > 1$. Thus, by the Dominated Convergence Theorem,
$$
\lim_{R \rightarrow +\infty} R \, \beta(R) = \frac{\lambda_\star}{4} \int_1^{+\infty} \frac{4}{z^2} \, dz = \lambda_\star,
$$
which concludes the proof of the proposition.
\ndim

Thanks to the various results displayed in the last subsections, the proof of Theorem~\ref{theorem1} is now complete.

\section{Proof of Theorem~\ref{theorem2}} \label{thm2sec}

\noindent

In this conclusive section, we finally address the proof of Theorem~\ref{theorem2}. Our argument essentially follows the lines of that displayed in~\cite[Section~5]{PSV13}. We stress that, aside from the obvious modifications due to the different framework in which our paper is set, we also correct some small mistakes present in~\cite{PSV13}.

\smallskip

Recalling definition~\eqref{u*def}, we have to prove that~$u^*$ is a class~A minimizer for~$\E_K$ and that it satisfies assertions~\eqref{extenests<1/2}-\eqref{extenvanishs>1/2}.

First of all, recall that~$u_0$ and, consequently,~$u^*$ are of class~$C^{1 + 2 s + \alpha}$, for some~$\alpha > 0$. Then, notice that
\begin{equation} \label{posparte}
\partial_{x_N} u^*(x) = \varpi u_0'(\varpi x_N) > 0 \quad \mbox{for any } x \in \R^N,
\end{equation}
and
\begin{equation} \label{limdire}
\lim_{x_N \rightarrow \pm \infty} u^*(x', x_N) = \lim_{x_N \rightarrow \pm \infty} u_0(\varpi x_N) = \pm 1 \quad \mbox{for any } x' \in \R^{N - 1}.
\end{equation}
Thus, by~\eqref{posparte},~\eqref{limdire} and Theorem~\ref{theorem3}, we are only left to show that~$u^*$ solves
\begin{equation} \label{ELND}
L_K u^* = W'(u^*) \quad \mbox{in } \R^N,
\end{equation}
to prove that~$u^*$ is a class~A minimizer for~$\E_K$. This is indeed quite straightforward. By substituting~$t := \varpi z_N$, we compute
\begin{equation} \label{Lu*Lu0}
\begin{aligned}
L_K u^*(x) & = \frac{1}{2} \int_{\R^N} \left( u^*(x + z) + u^*(x - z) - 2 u^*(x) \right) K(z) \, dz \\
%& = \frac{1}{2} \int_{\R} \left( u_0(\varpi (x_N - z_N)) + u_0(\varpi (x_N + z_N)) - 2 u_0(\varpi x_N) \right) \left( \int_{\R^{N - 1}} K(z', z_N) \, dz' \right) dz_N \\
& = \frac{1}{2 \varpi} \int_{\R} \left( u_0(\varpi x_N + t) + u_0(\varpi x_N - t) - 2 u_0(\varpi x_N) \right) \left[ \int_{\R^{N - 1}} K \left( z', \frac{t}{\varpi} \right) \, dz' \right] dt \\
& = \frac{1}{2} \int_{\R} \left( u_0( \varpi x_N + t) + u_0(\varpi x_N - t) - 2 u_0(\varpi x_N) \right) k(t) \, dt \\
& = L_k u_0(\varpi x_N),
\end{aligned}
\end{equation}
for any~$x \in \R^N$. Recall that the kernel~$k$ was defined in~\eqref{1d kernel}. Therefore, since~$u_0$ is a solution of
$$
L_k u_0 = W'(u_0) \quad \mbox{in } \R,
$$
by~\eqref{Lu*Lu0} we obtain
$$
L_K u^*(x) = L_k u_0(\varpi x_N) = W'(u_0(\varpi x_N)) = W'(u^*(x)) \quad \mbox{for any } x \in \R^N,
$$
which is~\eqref{ELND}.

%In the following we will almost always denote by~$x_e$ the projection~$e \cdot x$. Moreover, we will use the same letter~$C$ for various positive constants whose values are not important to us.
%
%First notice that, for any~$x \in \R^N$, with~$|x_e|$ suitably large, we have
%$$
%{\| \nabla u^* \|}_{L^\infty(B_{|x_e|/2}(x))} = {\| u_0' \|}_{L^\infty([x_e - |x_e|/2, x_e + |x_e|/2])} \le C |x_e|^{- 1 - 2 s}.
%$$
%Therefore, applying Lemma~\ref{enestlem} with~$\psi = u^*$ and~$\rho = |x_e| / 2$, we obtain
%\begin{align} \nonumber
%\int_{\R^N} {|u^*(x) - u^*(y)|}^2 K(x - y) dy & \le C \left( {\| \nabla u^* \|}_{L^\infty(B_{|x_e|/2}(x))}^2 |x_e|^{2 (1 - s)} + {\| u^* \|}_{L^\infty(\R^N)}^2 |x_e|^{- 2 s} \right) \\ \nonumber
%& \le C \left( |x_e|^{- 6 s} + |x_e|^{- 2 s} \right) \\ \label{1est}
%& \le C |x_e|^{- 2 s},
%\end{align}
%if~$|x_e|$ is large enough. On the other hand, if~$|x_e| \le R/2$, have
%\begin{align} \nonumber \label{2est}
%\int_{B_R^c} {|u^*(x) - u^*(y)|}^2 K(x - y) dy & \le \Lambda \int_{B_R^c} \frac{{|u^*(x) - u^*(y)|}^2}{|x - y|^{n + 2 s}} dy \\
%& \le \Lambda \int_{B_R^c} \frac{{|u_0(x_e) - u_0(y_e)|}^2}{|x_e - y_e|^{n + 2 s}} dy
%\end{align}

Thus, we are left to prove formulae~\eqref{extenests<1/2}-\eqref{extenvanishs>1/2}. In the remainder of the section, we will frequently denote with~$c$ any positive constant, whose value may change from line to line. Also, the radius~$R$ will be always implicitly assumed large.

Set
$$
I_{N, s}(R) := \int_{B_R} \int_{\R^N \setminus B_R} \frac{\left| u^*(x) - u^*(y) \right|^2}{|x - y|^{N + 2 s}} \, dx dy.
$$
First, we claim that
\begin{equation} \label{thm3claim1}
\mbox{if } s \in [1/2, 1), \mbox{ then } \lim_{R \rightarrow +\infty} \frac{I_{N, s}(R)}{R^{N - 1} \Psi_s(R)} = 0.
\end{equation}
and
\begin{equation} \label{thm3claim2}
\mbox{if } s \in (0, 1/2), \mbox{ then } I_{N, s}(R) \le c R^{N - 2 s}.
\end{equation}
Recall that~$\Psi_s$ was defined in~\eqref{Psidef}. Note that, thanks to the right-hand inequality in~\eqref{weakellipticity}, claim~\eqref{thm3claim1} would then imply formulae~\eqref{extenvanishs=1/2} and~\eqref{extenvanishs>1/2}, while~\eqref{thm3claim2} would yield~\eqref{extenests<1/2}.

To prove~\eqref{thm3claim1} and~\eqref{thm3claim2}, we write~$I_{N, s}(R) = S_{N, s}(R) + T_{N, s}(R)$, where
\begin{align*}
S_{N, s}(R) & := \int_{B_R} \int_{ \left( \R^N \setminus B_R \right) \cap \{ |x_N| \le R \} } \frac{|u^*(x) - u^*(y)|^2}{|x - y|^{N + 2 s}} \, dx dy, \\
T_{N, s}(R) & := \int_{B_R} \int_{\{ |x_N| > R \}} \frac{|u^*(x) - u^*(y)|^2}{|x - y|^{N + 2 s}} \, dx dy.
\end{align*}

First, we deal with the term~$T_{N, s}(R)$. We compute
\begin{align*}
T_{N, s}(R) & = \int_{-R}^R \int_{ \{ |x_N| > R \} } \int_{B_{\sqrt{R^2 - |y_N|^2}}'} \int_{\R^{N - 1}} \frac{|u^*(x) - u^*(y)|^2}{|x - y|^{N + 2 s}} \, dx' dy' dx_N dy_N \\
& = \int_{-R}^R \int_{ \{ |x_N| > R \} } \frac{|u_0(\varpi x_N) - u_0(\varpi y_N)|^2}{|x_N - y_N|^{N + 2 s}} \int_{B_{\sqrt{R^2 - |y_N|^2}}'} \left[ \int_{\R^{N - 1}} \frac{dx'}{\left( 1 + \frac{|x' - y'|^2}{|x_N - y_N|^2} \right)^{\frac{N + 2 s}{2}}} \right] dy' dx_N dy_N.
\end{align*}
If we change variables in the inner integral setting~$z' = (x' - y') / |x_N - y_N|$, we get
\begin{equation} \label{TNtech}
\begin{aligned}
T_{N, s}(R) & = \frac{\omega_{N - 1}}{\varpi^{2 s}} \int_{-R}^R \int_{ \{ |x_N| > R \} } \frac{|u_0(\varpi x_N) - u_0(\varpi y_N)|^2}{|x_N - y_N|^{1 + 2 s}} \left[ R^2 - |y_N|^2 \right]^{\frac{N - 1}{2}} dx_N dy_N \\
& \le \frac{\omega_{N - 1}}{\varpi^{2 s}} R^{N - 1} \int_{-R}^R \int_{ \{ |x_N| > R \} } \frac{|u_0(\varpi x_N) - u_0(\varpi y_N)|^2}{|x_N - y_N|^{1 + 2 s}} \, dx_N dy_N,
\end{aligned}
\end{equation}
recalling~\eqref{varpidef}. By exploiting~\eqref{N=1tailest}, we then get
\begin{equation} \label{TRest}
T_{N, s}(R) \le c R^{N - 2 s}.
\end{equation}

Now, we address the term~$S_{N, s}(R)$. We compute
\begin{align*}
S_{N, s}(R) & = \int_{-R}^R \int_{-R}^R \int_{B_{\sqrt{R^2 - |y_N|^2}}'} \int_{\R^{N - 1} \setminus B_{\sqrt{R^2 - |x_N|^2}}'} \frac{|u^*(x) - u^*(y)|^2}{|x - y|^{N + 2 s}} dx' dy' dx_N dy_N \\
& = \int_{-R}^R \int_{-R}^R \frac{|u_0(\varpi x_N) - u_0(\varpi y_N)|^2}{|x_N - y_N|^{N + 2 s}} \left[ \int_{B_{\sqrt{R^2 - |y_N|^2}}'} \int_{\R^{N - 1} \setminus B_{\sqrt{R^2 - |x_N|^2}}'} \frac{dx' dy'}{\left( 1 + \frac{|x' - y'|^2}{|x_N - y_N|^2} \right)^{\frac{N + 2 s}{2}}} \right] dx_N dy_N \\
& = \int_{-R}^R \int_{-R}^R \frac{|u_0(\varpi x_N) - u_0(\varpi y_N)|^2}{|x_N - y_N|^{2 - N + 2 s}} \, J_{\frac{N + 2 s}{2}, N - 1} \left( \frac{\sqrt{R^2 - |y_N|^2}}{|x_N - y_N|}, \frac{\sqrt{R^2 - |x_N|^2}}{|x_N - y_N|} \right) dx_N dy_N,
\end{align*}
where in the last line we changed variables by setting
$$
w' = \frac{x'}{|x_N - y_N|}, \quad z' = \frac{y'}{|x_N - y_N|},
$$
and the quantity~$J$ is defined in~\eqref{Jdef}. Applying then Corollary~\ref{Jcor}, we get\footnote{Observe that we use the estimate for~$J$ provided by Corollary~\ref{Jcor} with~$\alpha = (N + 2 s) / 2$ and~$N - 1$ in place of~$N$.}
\begin{equation} \label{SRdeco}
S_{N, s}(R) \le c_\delta \left( S_{N, s, \delta}^{(1)}(R) + S_{N, s}^{(2)}(R) \right),
\end{equation}
where
\begin{align*}
S_{N, s, \delta}^{(1)}(R) & := \int_{-R}^R \int_{-R}^R \frac{|u_0(\varpi x_N) - u_0(\varpi y_N)|^2}{|x_N - y_N|^{(1 - \delta) (N + 2 s)}} (R^2 - |y_N|^2)^{N - 1 - \delta \frac{N + 2 s}{2} } dx_N dy_N \\
S_{N, s}^{(2)}(R) & := \int_{-R}^R \int_{\{ |y_N| < |x_N| \}} \frac{|u_0(\varpi x_N) - u_0(\varpi y_N)|^2}{|x_N - y_N|^{1 + 2 s}} \left[ (R^2 - |y_N|^2)^{\frac{N - 1}{2}} - (R^2 - |x_N|^2)^{\frac{N - 1}{2}} \right] dx_N dy_N,
\end{align*}
the value~$\delta$ satisfies
\begin{equation} \label{thm3deltabounds}
\delta \in \left( \frac{N - 1}{N + 2 s}, \frac{N}{N + 2 s} \right),
\end{equation}
and~$c_\delta$ is a positive constant which may depend on~$N$,~$s$ and~$\delta$.
 
To estimate the first integral, we take
$$
\delta = \frac{2 N - 1}{2(N + 2 s)},
$$
which is clearly admissible for~\eqref{thm3deltabounds}. Accordingly, taking advantage of H\"{o}lder's inequality and the fact that~$|u_0| \le 1$,
\begin{align*}
S_{N, s, \frac{2 N - 1}{2 (N + 2 s)}}^{(1)}(R) & =\int_{-R}^R \int_{-R}^R \frac{|u_0(\varpi x_N) - u_0(\varpi y_N)|^2}{|x_N - y_N|^{\frac{1 + 4 s}{2}}} (R^2 - |y_N|^2)^{ \frac{2 N - 3}{4} } dx_N dy_N \\
& \le c R^{\frac{2 N - 3}{2}} \int_{-R}^R \int_{-R}^R \left( \frac{|u_0(\varpi x_N) - u_0(\varpi y_N)|^2}{|x_N - y_N|^{1 + 2 s}} \right)^{\frac{1 + 4 s}{2 (1 + 2 s)}} |u_0(\varpi x_N) - u_0(\varpi y_N)|^{\frac{1}{1 + 2 s}} dx_N dy_N \\
& \le c R^{\frac{2 N - 3}{2}} [u_0]_{H^s([-\varpi R, \varpi R])}^{\frac{1 + 4 s}{1 + 2 s}} \left( \int_{-R}^R \int_{-R}^R |u_0(\varpi x_N) - u_0(\varpi y_N)|^2 dx_N dy_N \right)^{\frac{1}{2 (1 + 2 s)}} \\
& \le c R^{\frac{2 N - 3}{2} + \frac{1}{1 + 2 s}} [u_0]_{H^s([- \varpi R, \varpi R])}^{\frac{1 + 4 s}{1 + 2 s}}.
\end{align*}
Recalling~\eqref{beta2est}, we compute
\begin{equation} \label{u0semiHsPsis}
[u_0]_{H^s([-\varpi R, \varpi R])}^2 \le \int_{-\varpi R}^{\varpi R} \left( \int_{\R} \frac{\left| u_0(r) - u_0(t) \right|^2}{|r - t|^{1 + 2 s}} \, dr \right) dt \le c \int_0^{\varpi R} \frac{dt}{1 + t^{2 s}} \le c \Psi_s(R).
\end{equation}
Accordingly,
\begin{equation} \label{SR1est}
S_{N, s, \frac{2 N - 1}{2 (N + 2 s)}}^{(1)}(R) \le c \begin{cases}
R^{N - 2 s} & \quad \mbox{if } s \in (0, 1/2) \\
R^{N - 1} \left( \log R \right)^{\frac{3}{4}} & \quad \mbox{if } s = 1/2 \\
R^{N - 1 - \frac{1}{2} \frac{2 s - 1}{1 + 2 s}} & \quad \mbox{if } s \in (1/2, 1).
\end{cases}
\end{equation}

The term~$S_{N, s}^{(2)}$ is more delicate. We start supposing~$N \ge 3$. Notice that, if~$0 \le a \le b$ and~$\beta \ge 1$, then
$$
b^\beta - a^\beta = \beta \int_a^b t^{\beta - 1} \, dt \le \beta b^{\beta - 1} (b - a).
$$
Applying this formula with~$\beta = (N - 1) / 2$, we get
\begin{align*}
(R^2 - |y_N|^2)^{\frac{N - 1}{2}} - (R^2 - |x_N|^2)^{\frac{N - 1}{2}} & \le \frac{N - 1}{2} (R^2 - |y_N|^2)^{\frac{N - 3}{2}} \left( |x_N|^2 - |y_N|^2 \right) \\
& \le (N - 1) R^{N - 2} |x_N - y_N|,
\end{align*}
if~$|y_N| \le |x_N| \le R$. Using the above estimate in combination with H\"older's inequality and~\eqref{u0semiHsPsis},
\begin{equation} \label{SR2estN>2}
S_{N, s}^{(2)}(R) \le c \begin{cases}
R^{N - 2 s} & \mbox{if } s \in (0, 1/2) \\
R^{N - 1} \sqrt{ \log R } & \mbox{if } s = 1/2 \\
R^{N - 1 - \frac{2 s - 1}{1 + 2 s}} & \mbox{if } s \in (1/2, 1).
\end{cases}
\end{equation}
We address the case~$N = 2$ in a slightly different way. First, fix any~$\mu \in (1, 2)$ and notice that, for any~$0 \le a \le b$,
\begin{align*}
\sqrt{b} - \sqrt{a} & = \frac{1}{2} \int_a^b \frac{dt}{\sqrt{t}} \le {\left( \int_a^b t^{- \frac{\mu}{2}} dt \right)}^{\frac{1}{\mu}} {\left( \int_a^b dt \right)}^{\frac{\mu - 1}{\mu}} \\
& = \frac{1}{2} {\left[ \frac{2}{2 - \mu} \left( b^{\frac{2 - \mu}{2}} - a^{\frac{2 - \mu}{2}} \right) \right]}^{\frac{1}{\mu}} {\left( b - a \right)}^{\frac{\mu - 1}{\mu}}.
%& \le C R^{\frac{2 - \mu}{2}} {\left( b - a \right)}^{\frac{\mu - 1}{\mu}}.
\end{align*}
Hence, by choosing e.g.~$\mu = 3/2$ we deduce that
\begin{align*}
\sqrt{R^2 - |y_2|^2}- \sqrt{R^2 - |x_2|^2} & \le c R^{\frac{2}{3}} {|x_2 - y_2|}^{\frac{1}{3}},
\end{align*}
and thus, arguing as for~\eqref{SR2estN>2},
\begin{equation} \label{SR2estN=2}
S_{2, s}^{(2)}(R) \le c \begin{cases}
R^{2 - 2 s} & \quad \mbox{if } s \in (0, 1/2) \\
R \left( \log R \right)^{\frac{5}{6}} & \quad \mbox{if } s = 1/2 \\
R^{1 - \frac{1}{3} \frac{2 s - 1}{1 + 2 s}} & \quad \mbox{if } s \in (1/2, 1).
\end{cases}
\end{equation}

By combining~\eqref{SR1est} and either~\eqref{SR2estN>2} or~\eqref{SR2estN=2}, by~\eqref{SRdeco} we conclude that
\begin{equation} \label{SRestsge1/2}
\lim_{R \rightarrow +\infty} \frac{S_{N, s}(R)}{R^{N - 1} \Psi_s(R)} = 0 \quad \mbox{if } s \in [1/2, 1),
\end{equation}
and
\begin{equation} \label{SRest<1/2}
S_{N, s}(R) \le c R^{N - 2 s} \quad \mbox{if } s \in (0, 1/2).
\end{equation}

Formulae~\eqref{TRest},~\eqref{SRestsge1/2} and~\eqref{SRest<1/2} imply claims~\eqref{thm3claim1} and~\eqref{thm3claim2}.

We now show that~\eqref{lowerextenests<1/2} is true. Recall that we prove its validity under the stronger assumption~\eqref{ellipticity} on~$K$. To check~\eqref{lowerextenests<1/2}, we use the identity displayed on the first line of~\eqref{TNtech} to write
$$
I_{N, s}(R) \ge T_{N, s}(R) = \frac{\omega_{N - 1}}{\varpi^{2 s}} \int_{-R}^R \int_{ \{ |x_N| > R \} } \frac{|u_0(\varpi x_N) - u_0(\varpi y_N)|^2}{|x_N - y_N|^{1 + 2 s}} \left[ R^2 - |y_N|^2 \right]^{\frac{N - 1}{2}} dx_N dy_N.
$$
By restricting the above integral to the values~$|x_N| \le R/2$ and recalling~\eqref{N=1tailtech}, we get
$$
I_{N, s}(R) \ge c R^{N - 1} \int_{-\frac{R}{2}}^{\frac{R}{2}} \int_{ \{ |x_N| > R \} } \frac{|u_0(\varpi x_N) - u_0(\varpi y_N)|^2}{|x_N - y_N|^{1 + 2 s}} \, dx_N dy_N \ge c R^{N - 2 s}.
$$
By~\eqref{ellipticity}, the left-hand inequality of~\eqref{extenests<1/2} then follows.

Finally, we head to the proof of~\eqref{totenconvs=1/2} and~\eqref{totenconvs>1/2}. Let now~$s \in [1/2, 1)$. Arguing as in~\eqref{Lu*Lu0} and changing variables appropriately, we get
\begin{align*}
\int_{\R^N} \int_{B_R} & {|u^*(x) - u^*(y)|}^2 K(x - y) \, dx dy \\
%& = \varpi \omega_{N - 1} R^{N - 1} \int_{-R}^R \left( \int_{\R} {|u_0(\varpi x_N) - u_0(\varpi y_N)|}^2 k \left( \varpi (x_N - y_N) \right) \, dy_N \right) {\left( 1 - \frac{|x_N|^2}{R^2} \right)}^{\frac{N - 1}{2}} dx_N \\
& = \frac{\omega_{N - 1} R^{N - 1}}{\varpi} \int_{-\varpi R}^{\varpi R} \left( \int_{\R} {|u_0(t) - u_0(r)|}^2 k \left( t - r \right) \, dr \right) {\left( 1 - \frac{t^2}{\varpi^2 R^2} \right)}^{\frac{N - 1}{2}} dt.
\end{align*}
Moreover, we easily compute
$$
\int_{B_R} W(u^*(x)) \, dx = \frac{\omega_{N - 1} R^{N - 1}}{\varpi} \int_{- \varpi R}^{\varpi R} W(u_0(t)) {\left( 1 - \frac{t^2}{\varpi^2 R^2} \right)}^{\frac{N - 1}{2}} dt.
$$
Hence, we write
\begin{equation} \label{totenconvtech}
\begin{aligned}
\frac{\E_K(u^*, B_R)}{R^{N - 1}} & = \frac{\omega_{N - 1}}{\varpi} \left[ \frac{1}{4} \int_{- \varpi R}^{\varpi R} \int_{\R} {|u_0(t) - u_0(r)|}^2 k(t - r) \, dr dt + \int_{- \varpi R}^{\varpi R} W(u_0(t)) \, dt \right] \\
& \quad + \theta_1(R) - \theta_2(R),
\end{aligned}
\end{equation}
where
\begin{align*}
\theta_1(R) & = \frac{1}{4 R^{N - 1}} \int_{\R^N \setminus B_R} \int_{B_R} {|u^*(x) - u^*(y)|}^2 K(x - y) \, dx dy \\
\theta_2(R) & = \frac{\omega_{N - 1}}{\varpi} \int_{- \varpi R}^{\varpi R} \alpha(t, \varpi R) \beta(t) \, dt,
\end{align*}
with
$$
\alpha(t, R') = 1 - {\left( 1 - \frac{t^2}{R^2} \right)}^{\frac{N - 1}{2}},
$$
and~$\beta$ as in~\eqref{betaest}, with~$k$ in place of~$K$.

Notice that
\begin{equation} \label{theta1van}
\lim_{R \rightarrow +\infty} \frac{\theta_1(R)}{\Psi_s(R)} = 0,
\end{equation}
by~\eqref{extenvanishs=1/2} or~\eqref{extenvanishs>1/2}. Furthermore, we claim that it also holds
\begin{equation} \label{theta2van}
\lim_{R \rightarrow +\infty} \frac{\theta_2(R)}{\Psi_s(R)} = 0.
\end{equation}
In order to check that~\eqref{theta2van} is valid, we distinguish between the two possibilities~$s = 1/2$ and~$s > 1/2$.

The latter case is easier. Indeed, when~$s > 1/2$, we know by~\eqref{betaest} that~$\beta \in L^1(\R)$. Since~$\alpha \le 1$, we may simply employ the Dominated Convergence Theorem to deduce~\eqref{theta2van}.

Conversely, when~$s = 1/2$ we need a more refined argument, inspired by~\cite[Lemma~4]{PSV13}. Write~$R' := \varpi R$. First we claim that, for any fixed~$\kappa \in (0, 1)$,
\begin{equation} \label{kappatech}
\lim_{R' \rightarrow +\infty} \frac{1}{\log R'} \int_{ \{ \kappa R' < |t| \le R' \} } \beta(t) \, dt = 0.
\end{equation}
Indeed, by~\eqref{betaest}, for any~$R' \ge 1$ we have
$$
\int_{ \{ \kappa R' < |t| \le R' \} } \beta(t) \, dt \le c_1 \int_{\kappa R'}^{R'} \frac{1}{1 + t} \, dt = c_1 \log \frac{1 + R'}{1 + \kappa R'} \le c_1 \log \frac{2}{\kappa},
$$
for some~$c_1 > 0$. From this,~\eqref{kappatech} clearly follows. In view of~\eqref{kappatech}, the way~$\alpha$ is defined and, again,~\eqref{betaest},
\begin{align*}
\lim_{R \rightarrow +\infty} \frac{\theta_2(R)}{\log R} & = \frac{\omega_{N - 1}}{\varpi} \lim_{R' \rightarrow +\infty} \frac{1}{\log R'} \left[ \int_{- \kappa R'}^{ \kappa R'} \alpha(t, R') \beta(t) \, dt + \int_{ \{ \kappa R' < |t| \le R' \} } \alpha(t, R') \beta(t) \, dt \right] \\
& \le c_2 \left[ 2 \lim_{R' \rightarrow +\infty} \frac{1}{\log R'} \left( \int_0^{ \kappa R'} \frac{dt}{1 + t} \right) \sup_{|\tau| \le \kappa R} \alpha(\tau, R') + \lim_{R' \rightarrow +\infty} \frac{1}{\log R'} \int_{ \{ \kappa R' < |t| \le R' \} } \beta(t) \, dt \right] \\
& = 2 c_2 \left[ 1 - \left( 1 - \kappa^2 \right)^{\frac{N - 1}{2}} \right],
\end{align*}
where~$c_2 > 0$ is independent of~$\kappa$. Since we may take~$\kappa$ as small as we like, we deduce that~\eqref{theta2van} is true also in this case.

By using~\eqref{theta1van} and~\eqref{theta2van} in~\eqref{totenconvtech}, it is easy to see that~\eqref{totenconvs=1/2} and~\eqref{totenconvs>1/2} are valid. Also,~\eqref{totenconvs=1/2b} follows from~\eqref{barEs=1/2} in Theorem~\ref{theorem1}, by noticing that if~$K$ satisfies~\eqref{Khom}, then the one dimensional kernel~$k$ defined by~\eqref{1d kernel} is of the type~\eqref{N=1Khom}, with~$\lambda_\star$ given by~\eqref{lambdastarthm3}. Indeed, using~\eqref{Khom} we compute
$$
k(t) = \frac{1}{\varpi} \int_{\R^{N - 1}} K \left( z', \frac{t}{\varpi} \right) \, dz' = \frac{\varpi^{N - 1 + 2 s}}{|t|^{N + 2 s}} \int_{\R^{N - 1}} K \left( \frac{\varpi z'}{t}, 1 \right) \, dz',
$$
for a.a.~$t \ne 0$. Changing now coordinates by setting~$y' = \varpi z' / t$, we get
$$
k(t) = \frac{\varpi^{2 s}}{|t|^{1 + 2 s}} \int_{\R^{N - 1}} K \left( y', 1 \right) \, dy' = \lambda_\star |t|^{- 1 - 2 s},
$$
and we are done.
This concludes the proof of Theorem~\ref{theorem2}.

\section*{Acknowledgements}

We would like to express our gratitude to Enrico Valdinoci for his precious advices and enduring support during the writing of this paper. 

%Bibliography ---------------------------------------------------------------------------------------------------------------------------------------------

\end{document}